\documentclass[leqno,12pt]{article}
\usepackage{latexsym}
\renewcommand{\epsilon}{\varepsilon}

\usepackage{amsmath}
\usepackage{amssymb}
\usepackage{ntheorem}
\usepackage[ansinew]{inputenc}
\usepackage{graphicx}
\usepackage{epsfig}
\usepackage{caption}

\usepackage[authoryear]{natbib}
\usepackage{amssymb}
\usepackage{setspace}
\usepackage{color}
\usepackage{amsmath}
\usepackage{graphicx}

\begin{document}
\newtheorem{theorem}{Theorem}[section]
\newtheorem{lemma}[theorem]{Lemma}
\newtheorem{rem}[theorem]{Remark}
\newtheorem{definition}[theorem]{Definition}
\newtheorem{assump}{Assumption}
\newtheorem{exam}[theorem]{Example}
\title{Additive  inverse regression models with convolution-type operators}

\author{Thimo Hildebrandt, Nicolai Bissantz, Holger Dette  \\
Ruhr-Universit\"at Bochum \\
Fakult\"at f\"ur Mathematik \\
44780 Bochum \\
Germany \\
{\small email: thimo.hildebrandt@ruhr-uni-bochum.de }\\
{\small  \qquad nicolai.bissantz@ruhr-uni-bochum.de }\\
{\small \quad holger.dette@ruhr-uni-bochum.de} \\
{\small FAX: +49 2 34 32 14 559}\\
}

 \maketitle

\begin{abstract}
\textcolor{white}{a}
\\
In a recent paper \cite{bb2008} considered the problem of nonparametric estimation in inverse regression models with convolution-type operators. For multivariate predictors nonparametric methods suffer from the curse of dimensionality and we consider inverse regression models with the additional qualitative assumption of additivity. In these models several additive estimators are studied. In particular, we  investigate estimators under the random design assumption which are applicable when observations are not available on a grid. Finally, we compare this estimator with the marginal integration and the non-additive estimator by means of a simulation study. It is demonstrated that the new method yields a substantial improvement of the currently available procedures.
\end{abstract}


Mathematical subject codes: primary, 62G08; secondary, 62G15, 62G20

\section{Introduction}
\def\theequation{1.\arabic{equation}}
\setcounter{equation}{0}

Inverse models have numerous applications in such important fields as biology, astronomy, economy or physics, where they have been
 intensively  studied  in a deterministic framework  [\cite{enghanneu1996}, \cite{saitoh1997}].
Recently inverse problems have also found considerable interest in the  statistical literature. These investigations
reflect  the demand in applications to quantify the uncertainty of estimates or to validate the model assumptions by
the construction of statistical confidence regions or  hypotheses tests, respectively [see  \cite{mairuy1996},  \cite{kaisom2005},
 \cite{bishohmunruy2007},   \cite{cavalier2008}, \cite{berbocdes2009}, \cite{berbocdes2009} or \cite{birbishol2010} among others].
In this paper we are interested in the convolution type inverse regression model
\begin{eqnarray}
Y &=& g(\textbf{z} )+\epsilon \label{model}
 =  \int_{\mathbb{R}^d} \psi(\textbf{z}-\textbf{t})\theta(\textbf{t}) d(\textbf{t}) +\epsilon
\end{eqnarray}
with a known function $\psi:\mathbb{R}^d\rightarrow \mathbb{R}$ [e.g. \cite{adorf1995}] and a centered noise term $\epsilon$. The goal of the experiment is to
recover  the signal $\theta: \mathbb{R}^d \rightarrow \mathbb{R}$ from  data $(\textbf{z}_{\textbf{1}},Y_{\textbf{1}}), \ldots , (\textbf{z}_{\textbf{n}},Y_{\textbf{n}}) $   which is closely
related to deconvolution [e.g. \cite{stef1990} and \cite{fan1991}]. Models of the type  (\ref{model}) have important
applications  in the recovery of images from  astronomical telescopes or fluorescence microscopes in biology. Therefore statistical inference for  the problem of estimating
the signal $\theta$ in model \eqref{model}  has become an important field of research in recent years, where the main focus is on a one dimensional predictor.  Bayesian methods have been
investigated in \cite{berbocdes2009} and \cite{kaisom2005} and  nonparametric methods have been  proposed by  \cite{mairuy1996}, \cite{cavalier2008} and \cite{bishohmunruy2007} among others.

In the present paper we investigate convergence properties of Fourier-based estimators for the function  $\theta$ with the following
purposes. Firstly, our research  is motivated by the fact that deconvolution problems often arise with a multivariate
predictor such as location and time. For this situation \cite{bb2008} proposed a nonparametric estimate of the signal $\theta$ and derived its asymptotic properties under rather strong assumptions. We will discuss the nonparametric estimation problem for the signal $\theta$
under substantially weaker assumptions. Secondly,
because nonparametric estimation usually suffers from the curse of dimensionality improved estimators incorporating qualitative assumptions such as  additivity or multiplicity are investigated under  the fixed and the random design assumption. While additive estimation has been intensively
discussed for direct  problems from different perspectives [see   \cite{linnie1995}, \cite{mamlinnie1999}, \cite{carhaemam2002},
\cite{henspe2005}, \cite{niespe2005}, \cite{dokkoo2000}, \cite{horlee2005}, \cite{leemampar2010}, \cite{detsch2011}]
- to our best knowledge - only one additive estimator is available for indirect inverse regression models so far where it is assumed that the observations are available on a grid [see \cite{birbishil2012}]. In this paper we are particularly interested in two alternative additive estimators. The first one is applicable if observations are available on a grid but has a substantially simpler structure than the method proposed by the last-named authors, which makes it very attractive for practitioners. Moreover, it also yields substantially more precise estimates than the method of \cite{birbishil2012}. The second estimator is additionally applicable in the case of random predictors.
\\Thirdly, we will also investigate the case
of correlated errors in the inverse regression model (\ref{model}), which has  - to our best knowledge - not been considered so far although it appears frequently in applications.
Finally, we do not assume that the kernel $\psi$ is periodic, which is a  common assertion in  inverse regression models with convolution operator [see e.g. \cite{cavtsy}]. Note that for many problems such as the reconstruction of  astronomical and biological images from telescopic and microscopic imaging devices
this assumption is unrealistic.

The remaining part of this paper is organized as follows.
In Section \ref{sec2} we introduce the necessary notation, different types of designs and estimators studied in this paper.
Section \ref{sec3} is devoted to the asymptotic properties of the estimators and we establish asymptotic normality
of all considered (appropriately standardized) statistics. In Section \ref{sec3a} we explain how the results are changing for dependent data while Section \ref{sec4} presents a small simulation study of the finite sample properties of the proposed methods. In particular we compare the new additive estimator with the currently available methods and demonstrate its superiority by a factor 6-8 with respect to mean squared error. Finally
all details regarding the proofs of our asymptotic results can be found in Section \ref{sec5}.

\section{Preliminaries} \label{sec2}
\def\theequation{2.\arabic{equation}}
\setcounter{equation}{0}

Recall the definition of model (\ref{model}) where we assume that the moments $E[\epsilon^k]$ exist for all $k \in \mathbb{N}$ such that $E[\epsilon]=0$ and $\sigma^2=E[\epsilon^2]>0$. For the sake of transparency we assume at this point that the errors corresponding to different predictors are independent - for the more general case of an error process with an MA($q$)-structure, see Section \ref{sec3a}.
We will investigate various estimators under two assumptions regarding
the explanatory variables $\textbf{z}$.
\begin{itemize}
\item[(FD)]  Under the  {\it fixed design} assumption we assume that observations are available on a grid of increasing size.
More precisely we consider a sequence $a_n \to 0$ as $n \to \infty$ and assume that  at each
location  $\textbf{z}_{\textbf{k}}=\frac{\textbf{k}}{na_n} \in  \mathbb{R}^d$ with $ \textbf{k}=(k_1,...,k_d) \in \{-n,...,n\}^d$ a pair of observations
$(\textbf{z}_{\textbf{k}},Y_{\textbf{k}}) $ is available in the model
\begin{eqnarray}
Y_{\textbf{k}}=g(\textbf{z}_{\textbf{k}})+\epsilon_{\textbf{k}} \label{model3}
= \int_{\mathbb{R}^d} \psi(\textbf{z}_{\textbf{k}}-\textbf{t})\theta(\textbf{t}) d\textbf{t} +\epsilon_{\textbf{k}} ,
\end{eqnarray}
where $\{\epsilon_\textbf{k} \: | \: \textbf{k} \in \{-n,...,n\}^d \}$ are independent and identically distributed random variables. Under this assumption the sample size is $N=(2n+1)^d$. Note that formally the random variables $\{Y_\textbf{k} \: | \: \textbf{k}\in \{-n,...,n\}^d \}$ form a triangular array, but we do not reflect this dependence in the notation. In other words we will use the notation $Y_\textbf{k}, \textbf{z}_\textbf{k}, \epsilon_\textbf{k}$ instead of $Y_{\textbf{k},n}, \textbf{z}_{\textbf{k},n}, \epsilon_{\textbf{k},n}$ throughout this paper.
\item[(RD)]  Under the  {\it random design} assumption we assume that the explanatory variables are realizations of
independent, identically distributed random variables
$\textbf{X}_{1,n},...,\textbf{X}_{n,n}$ with a density $f_n$. Again we will not reflect the triangular structure in the notation and use $Y_k , \textbf{X}_k,\epsilon_k$ and $f$ instead of $Y_{k,n}, \textbf{X}_{k,n}, \epsilon_{k,n}$ and $f_n$, respectively, that is
\begin{eqnarray}
Y_{k} = g(\textbf{X}_{k})+\epsilon_{k} \label{model2}
= \int_{\mathbb{R}^d} \psi(\textbf{X}_{k}-\textbf{t})\theta(\textbf{t}) d\textbf{t} +\epsilon_{k}; \quad k \in \{1,...,n\},
\end{eqnarray}
where $\epsilon_1,...,\epsilon_n$ are independent identically distributed random variables. Under this assumption the sample size is $N=n$.
\end{itemize}
We will use different estimators in both scenarios  \eqref{model3} and \eqref{model2}. Note that assumption (FD) assumes
that observations are available on a complete $d$-dimensional grid of length ${1\over na_n}$. In this case an estimator
of the signal $\theta$  has  also been studied by  \cite{bb2008}. The  estimator in model \eqref{model2}  under assumption (RD), which is proposed in the following section, could also be used if not all observations are
available on the grid.

\subsection{Unrestricted estimation for random design} \label{sec21}
Fourier-based estimators have been considered by numerous authors in the univariate case (e.g. \cite{digglehall93}, \cite{mairuy1996}, \cite{cavtsy} and \cite{bisdu2007}) and its generalization
to the multivariate case considered in the models \eqref{model3} and \eqref{model2} is straightforward.  For model \eqref{model3} a Fourier-based estimator is given by
\begin{eqnarray}
\hat{\theta}^{FD}(\textbf{x}^*) =\frac{1}{(2\pi)^d} \int_{\mathbb{R}^d} e^{-i\langle \textbf{w},\textbf{x}^*\rangle}\Phi_K(h\textbf{w}) \frac{\hat{\Phi}^{FD}(\textbf{w})}{\Phi_\psi(\textbf{w})} d\textbf{w}, \label{fest}
\end{eqnarray}
where
\begin{eqnarray*}
\hat{\Phi}^{FD}(\textbf{w}) &=& \frac{1}{n^da_n^d}\sum_{\textbf{k} \in \{-n,...,n\}^d} Y_\textbf{k} e^{i\langle \textbf{w},\textbf{z}_k\rangle }.
\end{eqnarray*}
denotes the empirical Fourier transform, $\langle \textbf{v}, \textbf{w}\rangle $ is the standard inner product of the vectors $\textbf{v},\textbf{w}\in \mathbb{R}^d$
and  $\Phi_K$  and $\Phi_\psi$   denote the Fourier transform of a kernel function $K$ and  the convolution function $\psi$ (which is assumed to be known),
respectively. Moreover,  in \eqref{fest} the quantity $h$ is a bandwidth converging to $0$ with increasing sample size.
\cite{birbishil2012} used this estimator to construct improved estimators under the qualitative assumption of additivity in the case of a
 fixed design. In Section \ref{sec22} we will propose an alternative additive estimator in the case of fixed design, which provides a notable improvement of the estimator proposed by the last named authors. \\
 For a random design we will use the same Fourier-based estimator as defined in  \eqref{fest},
 where the empirical Fourier transform $\hat{\Phi}^{FD}(w)$ in \eqref{fest} is replaced by
 \begin{eqnarray}
\hat{\Phi}^{RD}(\textbf{w})  &=& \frac{1}{n}\sum_{k =1}^n  e^{i\langle \textbf{w},\textbf{X}_k\rangle } \frac{Y_k}{\max \{f(\textbf{X}_k),f(\frac{1}{\textbf{a}_\textbf{n}}) \}}, \label{fourierrd}
\end{eqnarray}
$ 1/\textbf{a}_n=(1/a_n,...,1/a_n) \in \mathbb{R}^d $ and  $a_n$ is again a sequence  converging to $0$ with increasing sample size.
The resulting  estimator will be denoted by $\hat{\theta}^{RD}(\textbf{x}^*)$.
In \eqref{fourierrd} $f$ denotes the density of $\textbf{X}_1$
and we take the maximum of $f(\textbf{X}_k)$ and $f(\frac{1}{\textbf{a}_\textbf{n}})$ to ensure that the variance of $\hat{\theta}^{RD}(\textbf{x}^*)$  is bounded. We also note that the estimator $\hat{\theta}^{RD}$ admits the representation
\begin{eqnarray}
\hat{\theta}^{RD}(\textbf{x}^*) &=& \sum_{k=1}^n Y_k w_{n}(\textbf{x}^*,\textbf{X}_k), \label{linrep1}
\end{eqnarray}
where the weights are given by
\begin{eqnarray}
w_{n}(\textbf{x}^*,\textbf{X}_k) &=& \frac{1}{n\max\{f(\textbf{X}_k),f(\frac{1}{{\bf a_n}})\}(2\pi)^d} \int_{\mathbb{R}^d} e^{-i\langle \textbf{w},\textbf{x}^*-\textbf{X}_k \rangle} \frac{\Phi_K(h\textbf{w})}{\Phi_\psi(\textbf{w})} d\textbf{w}. \label{wRD}
\end{eqnarray}

\begin{rem}
{\rm Note that we use the same bandwidth for all components of the predictor. This assumption is made for the sake of a transparent presentation of the results. In applications the components of the vector $\textbf{x}$ represent different physical quantities such that different bandwidths have to be used. All results presented in this paper can be modified to this case with an additional amount of notation.}
\end{rem}

\subsection{Estimation of additive inverse regression models} \label{sec22}

It is well
known that in practical applications   nonparametric methods as introduced in Section \ref{sec21} suffer from the curse of dimensionality
and therefore do not yield precise estimates  of the signal $\theta$ with a multivariate predictor.
A   common approach in nonparametric statistics to deal with this problem is to postulate an additive
structure  of the signal  $\theta$, that is
\begin{eqnarray}
\theta(\textbf{x}^*)=\theta^{add}(\textbf{x}^*) := \theta_0^{add}+\sum_{j=1}^m \theta_{I_j}^{add}(\textbf{x}_{I_j}^*)  \label{thetaadd}
\end{eqnarray}
[see \cite{hastie1990}].
Here  $\{ I_1,...,I_m \} $ denotes a partition of the set  $ \{1,...,d\}$  with cardinalities $|I_j|=d_j$  and $\textbf{x}_{I_j}^*$ is the vector which includes all
components of the vector $\bf{x}^*$ with  corresponding indices $i \in I_j$. Furthermore
 $\theta_0^{add}$ is a constant and
 $\theta_{I_j}^{add}: \mathbb{R}^{d_j} \rightarrow \mathbb{R}$
 denote functions normalized such that  $$\int \theta_{I_j}^{add}(\textbf{x}) d(\textbf{x}) = 0  \quad  (j=1,...,m).$$
 Note that the  completely additive case is obtained for the choice $m=d$, that is $d_1=...=d_d=1$.
  In the case of direct regression models several
   estimation techniques such as marginal integration
    [see \cite{linnie1995}, \cite{carhaemam2002}, \cite{henspe2005}],  backfitting [\cite{mamlinnie1999}, \cite{niespe2005}]
have been proposed in the literature. Recently  the estimation problem of an additive (direct) regression model
has also found considerable interest in the context of quantile regression [see \cite{dokkoo2000}, \cite{goozer2003}, \cite{horlee2005},
\cite{leemampar2010},
\cite{detsch2011} among others] but - to our best knowledge - only one estimator has been proposed for additive inverse regression models under the assumption that observations are available on a grid [see \cite{birbishil2012}]. For this situation we will propose an alternative estimator in the following section, which yields an improvement by a factor 6-10 with respect to mean squared error (see our numerical results in Section \ref{sec4}).

To construct an estimator in the additive inverse regression model \eqref{thetaadd} with random design we apply the marginal integration method introduced in \cite{lin1995} to the statistic defined in \eqref{linrep1}. To this end we consider weighting functions $Q_{I_1},...,Q_{I_m}, Q_{I_j}:\mathbb{R}^{d_j}\rightarrow \mathbb{R}$ and define
\begin{eqnarray}
 Q(\textbf{x}^*) &=& Q_{I_1}(\textbf{x}_{I_1}^*)...Q_{I_m}(\textbf{x}_{I_m}^*) \nonumber
\\ Q_{I_j^c}(\textbf{x}_{I_j^c}^*) &=& Q_{I_1}(\textbf{x}_{I_1}^*)...Q_{I_{j-1}}(\textbf{x}_{I_{j-1}}^*)Q_{I_{j+1}}(\textbf{x}_{I_{j+1}}^*)...Q_{I_m}(\textbf{x}_{I_m}^*), \label{QIc}
\end{eqnarray}
where $I_j^c = \{1,\ldots , d\} \setminus {I_j}$.
With this notation we introduce the quantities
\begin{eqnarray}
\alpha_{j,Q_{I_j^c}}(\textbf{x}_{I_j}^*) &=& \int_{\mathbb{R}^{d-d_j}} \theta(\textbf{x}^*) dQ_{I_j^c}(\textbf{x}_{I_j^c}^*), \quad j=1,...,m, \label{a1}
\\ c &=& \int_{\mathbb{R}^d}\theta(\textbf{x}^*) dQ(\textbf{x}^*). \label{c}
\end{eqnarray}
Now let $\hat{\theta}^{RD}$ denote the unrestricted estimator introduced in Section 2.1 for the random design model, then the additive estimator for the signal $\theta$ is finally defined by
\begin{eqnarray}
\hat{\theta}^{add,RD}(\textbf{x}^*) &=&  \hat{\alpha}_{1,Q_{I_1^c}}(\textbf{x}_{I_1}^*)+...+\hat{\alpha}_{m,Q_{I_m^c}}(\textbf{x}_{I_m}^*)-(m-1)\hat{c} \label{margint}
\end{eqnarray}
where $\hat{c}$ and $\hat{\alpha}_{j,Q_{I_j^c}}$ denote estimates for the quantities $c$ and $\alpha_{j,Q_{I_j^c}}$ which are obtained by replacing in \eqref{a1} and \eqref{c} the signal $\theta$ by its estimator $\hat{\theta}^{full,RD}$, respectively. Recalling the definition of the unrestricted estimator in \eqref{fest} and \eqref{fourierrd}, we obtain from \eqref{a1} the representation
\begin{eqnarray}
\hat{\alpha}_{j,Q_{I_j^c}}(\textbf{x}_{I_j}^*) &=& \sum_{k =1}^n Y_k w_{n}^{add}(\textbf{x}_{I_j}^*,\textbf{X}_k), \label{thetaddrd}
\end{eqnarray}
where the weights are given by
\begin{eqnarray*}
w_{n}^{add}(\textbf{x}_{I_j}^*,\textbf{X}_k) &=& \frac{1}{nh^d(2\pi)^d} \int_{\mathbb{R}^d} e^{i\langle \textbf{w},\textbf{X}_k\rangle /h}e^{-i\langle \textbf{w}_{I_j},\textbf{x}_{I_j}^*\rangle /h}L_{I_j^c}\left(\frac{\textbf{w}_{I_j^c}}{h} \right) \frac{\Phi_K(\textbf{w})}{\Phi_\psi(\frac{\textbf{w}}{h})} d\textbf{w} \times\frac{1}{\max\{ f(\textbf{X}_k),f(\frac{1}{\textbf{a}_n})\}}
\end{eqnarray*}
and
\begin{eqnarray*}
L_{I_j^c}(\textbf{y}_{I_j^c}) &=& \int_{\mathbb{R}^{d-d_j}} e^{-i\langle \textbf{y}_{I_j^c},\textbf{x}_{I_j^c}^*\rangle }dQ_{I_j^c}(\textbf{x}_{I_j^c}^*).
\end{eqnarray*}

\subsection{An alternative additive estimator for a fixed design} \label{sec222}
In principle the marginal integration estimator could also be used under the fixed design assumption (FD)
and its asymptotic properties have been studied by \cite{birbishil2012}. However, it turns
out that for observations on a grid a simpler and more efficient estimator can be defined. This idea is closely related to the
backfitting approach. To be precise we note that the assumption of additivity for the signal $\theta$ implies additivity of the observable signal $g$ due to the linearity of the convolution operator. Hence, model \eqref{model3} is equivalent to
\begin{eqnarray}
\label{fixadd}
Y_{\textbf{k}} &=& g_0+g_{I_1}(\textbf{z}_{k_{I_1}})+...+g_{I_m}(\textbf{z}_{k_{I_m}})+\epsilon_{\textbf{k}},
\end{eqnarray}
where $g_0 = \int_{\mathbb{R}^d} \psi(\textbf{z}-\textbf{t}) \theta_0 d\textbf{t}$,
\begin{eqnarray}
\label{fixadd1}
g_{I_j}(\textbf{z}_{\textbf{k}_{I_j}}) =  \int_{\mathbb{R}^{d_j}} \psi_{I_j}(\textbf{z}_{\textbf{k}_{I_j}}-\textbf{t}_{I_j})\theta_{I_j}^{add}(\textbf{t}_{I_j}) d\textbf{t}_{I_j}
~~~~(j=1,\ldots ,m)
\end{eqnarray}
and $\psi_{I_1},...,\psi_{I_d}$ are the marginals  of $\psi$, that is
$$
\psi_{I_j}(\textbf{t}_{I_j}) ~=~\int_{\mathbb{R}^{d-d_j}} \psi(\textbf{t}) d\textbf{t}_{I_j^c}.
$$
Recall the definition of $\textbf{k}_{I_j}$ and $\textbf{k}_{I_j^c}$ as the $d_j$ and $(d-d_j)$-dimensional vector corresponding to the components $(k_l \:|\: l\in I_j)$ and $(k_l \:|\: l \in I_j^c)$ of the vector $\textbf{k}=(k_1,...,k_d)$, respectively. In order to define estimators of  these terms
 we consider the empirical Fourier transforms in dimension $d_j$
\begin{eqnarray*}
\hat{\Psi}_{{I_j}}(\textbf{w}) &=& \frac{1}{(na_n)^{d_j}} \sum_{\textbf{k}_{I_j} \in \{-n,...,n\}^{d_j}} Z_{\textbf{k}_{I_j}} e^{i \langle \textbf{w}, \textbf{z}_{\textbf{k}_{I_j}}\rangle }
~~(j=1,\ldots , m),
\end{eqnarray*}
where the random variables $Z_{\textbf{k}_{I_j}}$ are given by
\begin{eqnarray} \label{zmarg}
Z_{\textbf{k}_{I_j}} = \frac{1}{(2n+1)^{d-d_j}}\sum_{\textbf{k}_{I_j^c} \in \{-n,...,n\}^{d-d_j}} Y_{\textbf{k}} .
\end{eqnarray}
The additive estimator is now defined by
\begin{eqnarray}
\hat{\theta}^{add,FD}(\textbf{x}^*) &=&\hat{\theta}_0+ \hat{\theta}_{I_1}^{FD}(\textbf{x}_{I_1}^*)+...+\hat{\theta}_{I_m}^{FD}(\textbf{x}_{I_m}^*), \label{est}
\end{eqnarray}
where
\[ \hat{\theta}_0 =\frac{1}{n^d} \sum_{\textbf{k} \in \{-n,...,n\}^d} Y_{\textbf{k}} \]
\begin{eqnarray}  \label{estcomp}
\hat{\theta}_{I_j}^{FD}(\textbf{x}_{I_j}^*) &=& \frac{1}{(2\pi)^{d_j}} \int_{\mathbb{R}^{d_j}} e^{-i \langle \textbf{w}, \textbf{x}_{I_j}^*\rangle }
 \Phi_K(h\textbf{w}) \frac{ \hat{\Psi}_{g_{I_j}}(\textbf{w})}{\Phi_{\psi_{I_j}}(\textbf{w})} d\textbf{w}
 ~~~~ (j=1,\ldots , m).
\end{eqnarray}
 Note that by the lattice structure the statistic $Z_{\textbf{k}_{I_j}}$ in \eqref{zmarg}
is a $\sqrt{n^{d-d_j}}$-consistent estimator of $g_{I_j}(\textbf{z}_{\textbf{k}_{I_j}})$. Therefore the deconvolution
problem for the $j$-th component is reduced to a problem in  dimension $d_j$ and  the estimator $\hat{\theta}_{I_j}^{FD}(\textbf{x}_{I_j}^*)$ can be rewritten as
\begin{eqnarray}
\hat{\theta}_{I_j}^{FD}(\textbf{x}_{I_j}^*) &=& \sum_{\textbf{k}_{I_j} \in \{-n,...,n\}^{d_j}} Z_{\textbf{k}_{I_j}} w_{\textbf{k}_{I_j},n}(\textbf{x}_{I_j}^*), \label{sumFD}
\end{eqnarray}
where the  weights $w_{\textbf{k}_{I_j},n}$ are defined by
\begin{eqnarray}
w_{\textbf{k}_{I_j},n}(\textbf{x}_{I_j}^*) = \frac{1}{(nha_n2\pi)^{d_j}} \int_{\mathbb{R}^{d_j}} e^{-i\langle \textbf{w}, (\textbf{x}_{I_j}^*-\textbf{z}_{\textbf{k}_{I_j}})\rangle /h}  \frac{\Phi_K(\textbf{w}) }{\Phi_{\psi_{I_j}}(\frac{\textbf{w}}{h})} d\textbf{w}. \label{wfdadd}
\end{eqnarray}

\subsection{Technical Assumptions}
In the following Section we will derive important asymptotic properties of the proposed estimators. For this purpose the following assumptions are required, where different statements in the following discussion require different parts of these assumptions. Throughout this paper $\parallel . \parallel$ denotes the Euclidean norm and the symbol $a_n \sim b_n$ means that $\lim_{n \rightarrow \infty} a_n/b_n =c$ for some positive constant $c$.

\begin{assump}
\label{asspsiadd}
\color{white}{a} \color{black}
\begin{enumerate}
\item[{\rm (A)}] {\rm Under the random design assumption the Fourier transform $\Phi_\psi$  of the function $\psi$ satisfies (as $h\rightarrow 0$)
\begin{eqnarray*}
\int_{\mathbb{R}^d} \frac{|\Phi_K(\textbf{w})|}{|\Phi_{\psi}(\frac{\textbf{w}}{h})|} d\textbf{w} &\le& C_1h^{-\beta} ~,~~
 \int_{\mathbb{R}^d} \frac{|\Phi_K(\textbf{w})|^2}{|\Phi_{\psi}(\frac{\textbf{w}}{h})|^2} d\textbf{w}  \sim  C_2h^{-2\beta}
\end{eqnarray*}
for some $\beta > 0$ and constants $C_1,C_2 > 0$.}
\item[{\rm (B)}] {\rm Under the fixed design and additivity assumption the Fourier transforms $\Phi_{\psi_{I_j}}$ of the marginals $\psi_{I_j}$ of $\psi$  satisfy
\begin{eqnarray*}
\int_{\mathbb{R}^{d_j}} \frac{|\Phi_K(\textbf{w})|}{|\Phi_{\psi_{I_j}}(\frac{\textbf{w}}{h})|} d\textbf{w} &\le& C_1h^{-\beta_j} ~,~~
 \int_{\mathbb{R}^{d_j}} \frac{|\Phi_K(\textbf{w})|^2}{|\Phi_{\psi_{I_j}}(\frac{\textbf{w}}{h})|^2} d\textbf{w} \sim C_2h^{-2\beta_j}
\end{eqnarray*}
for some $\beta_j >0$ $(j=1,\dots,m)$ and constants  $C_1,C_2 >0$.}
\end{enumerate}
\end{assump}

\begin{assump}
\label{assk}
\textcolor{white}{a}
\begin{enumerate}
\item[{\rm (A)}] {\rm Under the random design assumption the Fourier transform $\Phi_K$ of the kernel $K$  in \eqref{fest} is symmetric, supported on the cube $[-1,1]^d$ and there exists a constant $b \in (0,1]$ such that $\Phi_K(\textbf{w})=1$ for $\textbf{w} \in [-b,b]^d,b> 0,$ and $|\Phi_K(\textbf{w})| \le 1$ for all $\textbf{w} \in [-1,1]^d.$}
\item[{\rm (B)}] {\rm Under the fixed design and additivity assumption the Fourier transform $\Phi_K$ of the kernel $K$ is symmetric and supported on $[-1,1]^{d_j}$ and there exists a constant $b \in (0,1]$ such that $\Phi_K(\textbf{w})=1$ for $\textbf{w} \in [-b,b]^{d_j},b> 0,$ and $|\Phi_K(\textbf{w})| \le 1$ for all $\textbf{w} \in [-1,1]^{d_j}$ for all $j=1,...,m$.}
\end{enumerate}
\end{assump}

\begin{assump}
\label{asstheta}
\color{white}{a} \color{black}
\begin{enumerate}
\item[{\rm (A)}] {\rm The Fourier transform $\Phi_{\theta}$ of the signal $\theta$ in model \eqref{model} exists and satisfies
\begin{eqnarray*}
\int_{\mathbb{R}^d} |\Phi_{\theta}(\textbf{w})| \parallel \textbf{w} \parallel^{s-1} d\textbf{w} < \infty \quad \mbox{for some} \quad s >1.
\end{eqnarray*}}
\item[{ \rm (B)}] {\rm  The function $g$ in model \eqref{model} satisfies
\begin{eqnarray*}
\int_{\mathbb{R}^d} |g(\textbf{z})| \parallel \textbf{z} \parallel^r d\textbf{z} < \infty
\end{eqnarray*}
for some $r>0$ such that $a_n^r = O(h^{\beta+d+s-1})$.}
\item[{\rm (C)}] {\rm The Fourier transforms $\Phi_{\theta_{I_1}^{add}},...,\Phi_{\theta_{I_m}^{add}}$ of the functions $\theta_{I_1}^{add},...,\theta_{I_m}^{add}$ in the additive model \eqref{thetaadd}  satisfy
\begin{eqnarray*}
\int_{\mathbb{R}^d} |\Phi_{\theta_{I_j}^{add}}(\textbf{w})| \parallel \textbf{w} \parallel^{s-1} d\textbf{w} < \infty \quad \mbox{for some} \quad s>1 \mbox{ and } j=1,...,m.
\end{eqnarray*}}
\item[{\rm (D)}] {\rm The functions $g_{I_1},...,g_{I_m}$ defined in \eqref{fixadd1} satisfy
\begin{eqnarray*}
\int_{\mathbb{R}^{d_j}} |g_{I_j}(\textbf{z})| \parallel \textbf{z} \parallel^r d\textbf{z} > \infty \quad \mbox{for} \quad j=1,...,m
\end{eqnarray*}
for some $r>0$ such that  $a_n^{r-d_j} = O(h^{\beta_j+s+d_j-1})$.}
\end{enumerate}
\end{assump}

\begin{assump}
\label{assrv}
 {\rm For each $n\in \mathbb{N}$ let $\textbf{X}_1,...,\textbf{X}_n$ denote independent identically distributed d-dimensional  random variables with density $f$ (which may depend on $n$) such that $f(\textbf{x}) \not= 0$ for all $\textbf{x}\in [-\frac{1}{a_n},\frac{1}{a_n}]^d$. Furthermore we assume, that for sufficiently large $n \in \mathbb{N}$
\[f(\textbf{x}) \ge f(\frac{1}{\bf{a}_n}) \quad \mbox{ for } \quad \textbf{x} \in [-\frac{1}{a_n},\frac{1}{a_n}]^d.\] }
\end{assump}
The final assumption is required for the marginal integration estimator and is an extension of Assumption \ref{asspsiadd}. For a precise statement we define for $\textbf{y} \in \mathbb{R}^{d-d_j}$
\begin{eqnarray}
L_{I_j^c}(\textbf{y}) = \int_{\mathbb{R}^{d-d_j}} e^{-i\langle \textbf{y},\textbf{x}_{I_j^c}\rangle } dQ_{I_j^c}(\textbf{x}_{I_j^c}) \label{LI}
\end{eqnarray}
where $Q_{I_j^c}(\textbf{x}_{I_j^c})$ as defined in \eqref{QIc}.
\\
\begin{assump}
\label{asspsi}
{\rm There exist positive constants $\gamma_1,...,\gamma_m$ such that the Fourier transform $\Phi_{\psi}$ of the convolution function $\psi$ satisfies
\begin{enumerate}
\item[{ \rm (A)}] $\int_{\mathbb{R}^d}\left|L_{I_j^c}\left(\frac{\textbf{w}_{I_j^c}}{h}\right)\right|^2 \frac{|\Phi_K(\textbf{w})|^2}{|\Phi_{\psi}(\frac{\textbf{w}}{h})|^2} d\textbf{w} \sim C_3h^{-2\beta+\gamma_j} \quad (j=1,...,m)$
\item[{\rm (B)}] $ \int_{\mathbb{R}^d} \left| \sum_{j=1}^m e^{-i\langle \textbf{w}_{I_j},\textbf{x}_{I_j}^*\rangle /h}L_{I_j^c}\left(\frac{\textbf{w}_{I_j^c}}{h} \right)\right|^2 \frac{|\Phi_K(\textbf{w})|^2}{|\Phi_\psi(\frac{\textbf{w}}{h})|^2} d\textbf{w} \sim C_4 h^{-2\beta+\gamma_{min}}, \mbox{ where } \gamma_{min}=\min_{j=1}^m \gamma_j$
\item[{\rm (C)}] $\int_{\mathbb{R}^d} \big(\prod_{j=1}^m \big|L_{I_j^c}\big( \frac{\textbf{w}_{I_j^c}}{h}\big)\big|^2\big)\frac{|\Phi_K(\textbf{w})|^2}{|\Phi_\psi(\frac{\textbf{w}}{h})|^2} d\textbf{w} = o\left(h^{-2\beta+\gamma_{min}}\right).$
\end{enumerate}
 }
\end{assump}

\begin{rem}
\color{white}{a} \color{black}
\begin{enumerate}
\item  {\rm The common assumption on the convolution function $\psi$ is
\begin{eqnarray}
\Phi_\psi(\textbf{w}) \parallel \textbf{w} \parallel^\beta \rightarrow C \quad \textbf{w} \rightarrow \infty, \label{weakbb}
\end{eqnarray}
[see \cite{bb2008}]. Assumption 1 is substantially weaker because we do not assume $\Phi_\psi$ to be  asymptotically radial-symmetric. It is satisfied for many commonly used convolution functions such as the multivariate Laplace density, the density of several Gamma distributions such as the Exponential distribution for which \eqref{weakbb} does not hold.
\item Assumptions \ref{asstheta}(A) and \ref{asstheta}(B) will not be required for the new additive estimator introduced in Section 2.2 under the fixed design assumption. As a consequence the asymptotic theory for the new estimator in the completely additive case $m=d$ $(d_1=...=d_m=1)$ does not require the additive functions  to have compact support as it is assumed in \cite{birbishil2012}.
\item Assumptions \ref{asstheta}(B) and \ref{asstheta}(D) are needed for the computation of the bias, where we have to ensure that $g(\textbf{x})$ converges sufficiently fast to zero as $x\rightarrow \infty$. Note that we only observe data on the cube $[ -\frac{1}{a_n},\frac{1}{a_n} ]^d$.
\item
The results of this Section can be extended to multiplicative signals of the form
 \begin{eqnarray}
\theta^{mult}(\textbf{x}^*) &=& \prod_{j=1}^m \theta_{I_j}^{mult}(\textbf{x}_{I_j}^*).
\end{eqnarray}
The details are omitted for the sake of brevity.
}
\end{enumerate}
\end{rem}

\begin{exam}
{ \rm In order to demonstrate that the assumptions are satisfied in many cases of practical importance we consider exemplarily Assumptions \ref{asspsiadd} and \ref{asspsi} and a two dimensional additive signal that is $\textbf{x}=(x_1,x_2)$,
\[ \theta(x_1,x_2)=\theta_1(x_1)+\theta_2(x_2), \]
$(I_1=I_2^c=\{1\},I_2=I_1^c=\{2\})$. For the convolution function in \eqref{model} and the weight \eqref{QIc} we choose
\allowdisplaybreaks{
\begin{eqnarray*}
 \psi(\textbf{x}) &=& \frac{\lambda^2}{4}e^{-\lambda(|x_1|+|x_2|)}
\\ Q(\textbf{x}) &=& 1_{[-1,1]^2}(\textbf{x}),
\end{eqnarray*}
respectively, and the kernel $K$ is given by
\[K(\textbf{x}) = \frac{\sin(x_1)\sin(x_2)}{\pi^2x_1x_2}.\]
 The integrals in Assumptions \ref{asspsiadd} and \ref{asspsi} are therefore obtained by a straightforward calculation \footnotesize
 \begin{eqnarray*}
\int_{\mathbb{R}^2} \frac{|\Phi_K(\textbf{w})|}{|\Phi_{\psi}(\frac{\textbf{w}}{h})|} d\textbf{w} &=& \int_{[-1,1]^2} \left(1+\frac{w_1^2}{h^2}\right)\left(1+\frac{w_2^2}{h^2}\right) d\textbf{w} = \left(\frac{2}{3h^2}+2\right)^2
\end{eqnarray*}
 \begin{eqnarray*}
 \int_{\mathbb{R}^2} \frac{|\Phi_K(\textbf{w})|^2}{|\Phi_{\psi}(\frac{\textbf{w}}{h})|^2} d\textbf{w}&=& \int_{[-1,1]^2} \left(1+\frac{w_1^2}{h^2}\right)^2\left(1+\frac{w_2^2}{h^2}\right)^2 d\textbf{w}= \left( \frac{2}{5h^4}+\frac{4}{3h^2}+2\right)^2
\end{eqnarray*}
 \begin{eqnarray*}
 \int_{\mathbb{R}^2}\left|L_{1}\left(\frac{w_{1}}{h}\right)\right|^2 \frac{|\Phi_K(\textbf{w})|^2}{|\Phi_{\psi}(\frac{\textbf{w}}{h})|^2} d\textbf{w} &=& \int_{[-1,1]^2} \frac{4h^2|\sin\left(\frac{w_1}{h}\right)|^2\left(1+\frac{w_1^2}{h^2}\right)^2\left(1+\frac{w_2^2}{h^2}\right)^2}{w_1^2} d\textbf{w}
= \frac{8}{15h^6}+o(h^{-6}),
\end{eqnarray*}
\normalsize
where we define $\mbox{Si}(x) =\int_0^x \frac{\sin(t)}{t} dt$.
\footnotesize
 \begin{eqnarray*}
 \int_{\mathbb{R}^2} \left| \sum_{j=1}^2 e^{-i\langle w_{I_j},\textbf{x}_{I_j}^*\rangle /h}L_{I_j^c}\left(\frac{w_{I_j^c}}{h} \right)\right|^2 \frac{|\Phi_K(\textbf{w})|^2}{|\Phi_\psi(\frac{\textbf{w}}{h})|^2} d\textbf{w} &=&h^2\int_{[-1,1]^2}  \left| e^{-iw_1x_1/h}\frac{\sin\left(\frac{w_2}{h}\right)}{w_2}+ e^{-iw_2x_2/h}\frac{\sin\left(\frac{w_1}{h}\right)}{w_1}\right|^2
\\ && \times \left(1+\frac{w_1^2}{h^2}\right)^2\left(1+\frac{w_2^2}{h^2}\right)^2 d\textbf{w}
 =  \frac{16}{15h^6} + o\left(h^{-6}\right)
\end{eqnarray*}
 \begin{eqnarray*}
 \int_{\mathbb{R}^2}\prod_{j=1}^2\left|L_{I_j}\left(\frac{w_{I_j}}{h}\right)\right|^2 \frac{|\Phi_K(\textbf{w})|^2}{|\Phi_{\psi}(\frac{\textbf{w}}{h})|^2} d\textbf{w} &=& \Big(\int_{[-1,1]} \frac{4h^2|\sin\left(\frac{w_1}{h}\right)|^2\left(1+\frac{w_1^2}{h^2}\right)^2}{w_1^2} dw_1\Big)^2
= \frac{16}{9h^4}+o(h^{-4})
\end{eqnarray*}
}}
\end{exam}
\normalsize

\section{Asymptotic properties} \label{sec3}
\def\theequation{3.\arabic{equation}}
\setcounter{equation}{0}

\subsection{Unrestricted estimator}
In the following we discuss the weak convergence of the unrestricted estimator $\hat{\theta}^{RD}$ for the signal $\theta$. In the case of a fixed design on a grid (assumption (FD)) the asymptotic properties of this estimator have been studied in \cite{bb2008}. Therefore we restrict ourselves to model \eqref{model2}  corresponding to the random design assumption, for which the situation is substantially more complicated. Here the estimator is given by
\begin{eqnarray}
\\
 \hat{\theta}^{RD}(\textbf{x}^*) =\frac{1}{nh^d(2\pi)^d} \sum_{k =1}^n  \int_{\mathbb{R}^d} e^{-i\langle \textbf{w},\textbf{x}^*-\textbf{X}_k\rangle /h}\frac{\Phi_K(\textbf{w})}{\Phi_\psi(\frac{\textbf{w}}{h})} d\textbf{w} \frac{Y_k}{\max\{ f(\textbf{X}_k),f(\frac{1}{\textbf{a}_\textbf{n}})\}} \nonumber  \label{thetfullrd}
\end{eqnarray}
and its asymptotic properties are described in our first main result which is proved in the appendix. Throughout this paper the symbol $\Rightarrow$ denotes weak convergence.

\begin{theorem}
\label{theoaddrd}
Consider the inverse regression model \eqref{model2} under the random design assumption (RD). Let Assumptions {\rm \ref{asspsiadd}(A), \ref{assk}, \ref{asstheta}(A), \ref{asstheta}(B), \ref{assrv} \textit{and} \ref{asspsi}} be fulfilled and $h \rightarrow 0$ and $a_n\rightarrow 0$ as $n \rightarrow \infty$ such that
\[ n^{1/2}h^{\beta+d/2}f(\textbf{a}_\textbf{n}^{-1})^{1/2} \rightarrow  \infty \qquad\mbox{ and } \qquad n^{1/2}h^{3d/2}f(\textbf{a}_\textbf{n}^{-1})^{3/2} \rightarrow \infty.\]
 Furthermore, assume that the errors in model \eqref{model2} are independent, identically distributed with mean zero and variance $\sigma^2$. Then
\begin{eqnarray}
 V_1 ^{-1/2} \big(\hat{\theta}^{RD}(\textbf{x}^*)-E[\hat{\theta}^{RD}(\textbf{x}^*)]\big) \Rightarrow {\cal{N}}(0,1), \label{weak1}
\end{eqnarray}
where $E[\hat{\theta}^{RD}(\textbf{x}^*)]=\theta(\textbf{x}^*)+O(h^{s-1})$ and the normalizing sequence
\begin{eqnarray} \label{v1}
 V_1 =\frac{1}{n(2\pi)^{2d}}\int_{\mathbb{R}^d}\left( \int_{\mathbb{R}^d} e^{-i\langle \textbf{s},(\textbf{x}^*-\textbf{y})\rangle }\frac{\Phi_K(\textbf{hs}) }{\Phi_{\psi}(\textbf{s})} d\textbf{s} \right)^2\frac{(\sigma^2+g^2(\textbf{y}))f(\textbf{y})}{\max\big \{f(\textbf{\textbf{y}}),f(\frac{1}{\textbf{a}_\textbf{n}})\big\}^2}d\textbf{y}
 \end{eqnarray}
is bounded by
\begin{eqnarray}
C_ln^{1/2}h^{d/2+\beta}f(\textbf{a}_\textbf{n}^{-1})^{1/2}\le V_1^{-1/2}  \le C_un^{1/2}h^{d/2+\beta}. \label{convrate}
\end{eqnarray}
\end{theorem}

\begin{rem} \label{rem3.2}
{ \rm
Note that the rate of convergence in Theorem \ref{theoaddrd} depends sensitively on the design density. We demonstrate this by providing two examples, one for the fastest and one for the slowest possible rate. First, assume that the predictors are uniformly distributed on the cube $[-\frac{1}{a_n},\frac{1}{a_n}]^d$ and that the convolution function is the $d$-dimensional Laplace density function. This yields $\beta=2d$ in Assumption 1 and we get a rate of convergence of order $n^{1/2}h^{5d/2} a_n^{d/2}$, which is exactly the lower bound in Theorem \ref{theoaddrd} and coincides with the rate in the fixed design case. However, a rate of order $n^{1/2}h^{5d/2}$ is obtained for the design density
\begin{eqnarray*}
f(x_1,...,x_d)&=& \prod_{k=1}^d g_{a,b}(x_k),
\end{eqnarray*}
where the function $g_{a,b}:\mathbb{R}\rightarrow \mathbb{R}$  is defined by
\begin{eqnarray*}
g_{a,b}(x)=\begin{cases}
 a,  & \text{if } x\in[-1,1]\\
  \frac{a}{| x |^b}, & \text{else },
\end{cases}
\end{eqnarray*}
and the parameters $a$ and $b$ are given by $b>1,a=(2+\frac{2}{b-1})^{-1}$. In this case we have
\[V_1^{-1/2} \sim n^{-1/2}h^{-5d/2}+n^{-1/2}h^{-2d}a_n^{(-b+1)/2}.\]
For the choice $h=o(a_n^{b-1})$ we therefore obtain $V_1^{-1/2} \sim n^{-1/2}h^{-5d/2}$. }
\end{rem}

\subsection{Additive estimation for random design}

In this Section we consider the marginal integration estimator $\hat{\theta}^{add,RD}$ defined in \eqref{margint} under the random design assumption. Lemma \ref{lemaddrd} below gives the asymptotic behaviour of the $j$-th component $\hat{\alpha}_{j,Q_{I_j^c}}$ and Theorem \ref{theoremaddrd} the asymptotic distribution of $\hat{\theta}^{add,RD}$. The proofs are complicated and also deferred to Section \ref{sec5}.

\begin{lemma}
\label{lemaddrd}
If Assumptions {\rm \ref{asspsiadd}(A), \ref{assk}, \ref{asstheta}(C), \ref{asstheta}(D), \ref{assrv} \textit{and} \ref{asspsi}} are satisfied and \[n^{1/2}h^{\beta+d/2-\gamma_j/2}f(\textbf{a}_\textbf{n}^{-1})^{1/2} \rightarrow \infty \qquad \mbox{ and } \qquad n^{1/2}h^{3/2(d-\gamma_j)}f(\textbf{a}_\textbf{n}^{-1}) \rightarrow \infty\]
 as $n \rightarrow \infty$. Then the appropriately standardized estimator $\hat{\alpha}_{j,Q_{I_j^c}}(\textbf{x}_{I_j}^*)$ defined in \eqref{thetaddrd} converges weakly to a standard normal distribution, that is
\begin{eqnarray}
 V_2^{-1/2}  \big(\hat{\alpha}_{j,Q_{I_j^c}}(\textbf{x}_{I_j}^*)-E[\hat{\alpha}_{j,Q_{I_j^c}}(\textbf{x}_{I_j}^*)])\big)
 \Rightarrow {\cal{N}}(0,1) \label{weakalpha}
\end{eqnarray}
for $j=1,...,m$, where $E[\hat{\alpha}_{j,Q_{I_j^c}}(\textbf{x}_{I_j}^*)]=\alpha_{j,Q_{I_j^c}}(\textbf{x}_{I_j}^*)+O(h^{s-1})$ and the standardizing factor
\begin{eqnarray*}
V_2 &=& \frac{1}{n(2\pi)^d}\int_{\mathbb{R}^d} \left(\int_{\mathbb{R}^d} e^{-i\langle \textbf{w},\textbf{x}\rangle } e^{i\langle \textbf{w}_{I_j},\textbf{x}_{I_j}\rangle }L_{I_j^c}\big(\textbf{w}_{I_j^c}\big)\frac{\Phi_K(h\textbf{w})}{\Phi_\psi(\textbf{w})} d\textbf{w}\right)^2 \frac{(\sigma^2+g(\textbf{x})^2)f(\textbf{x})}{\max\{f(\textbf{x}),f(\frac{1}{\textbf{a}_\textbf{n}})\}^2} d\textbf{x}.
\end{eqnarray*}
satisfies
\begin{eqnarray*}
C_ln^{1/2}h^{d/2+\beta-\gamma_j}f(\textbf{a}_\textbf{n}^{-1})^{1/2}\le V_2^{-1/2}  \le C_u n^{1/2}h^{d/2+\beta-\gamma_j}.
\end{eqnarray*}
\end{lemma}

\begin{rem}
{\rm
Similar to the unrestricted case, the rate of convergence depends on the design density $f$. Note that under the given assumptions the rate of convergence of the  estimator $\hat{\alpha}_{j,Q_{I_j^c}}$ is by the factor $h^{\gamma_j}$ faster than the rate of the unrestricted estimator.}
\end{rem}

\begin{theorem}
\label{theoremaddrd}
If Assumptions {\rm \ref{asspsiadd}(A), \ref{assk}, \ref{asstheta}(C), \ref{asstheta}(D), \ref{assrv} \textit{and} \ref{asspsi}} are satisfied and \begin{eqnarray*}
&&nh^{\beta+(3d+\gamma_{min}/)2}f(\textbf{a}_\textbf{n}^{-1})^2 \rightarrow \infty,  \quad n^{1/2}h^{\beta+(d-\gamma_{min})/2}f(\textbf{a}_\textbf{n}^{-1})^{1/2} \rightarrow \infty
\\&& n^{1/2}h^{3/2(d-\gamma_j)}f(\textbf{a}_\textbf{n}^{-1})^3 \rightarrow \infty \quad (j=1,...,m)
\end{eqnarray*}
 as $n \rightarrow \infty$, then the appropriately standardized additive estimator $\hat{\theta}^{add,RD}$ converges weakly to a standard normal distribution, that is
\begin{eqnarray}
V_3^{-1/2} \big( \hat{\theta}^{add,RD}(\textbf{x}^*)-E[\hat{\theta}^{add,RD}(\textbf{x}^*)] \big) \Rightarrow {\cal{N}}(0,1), \label{weak2}
\end{eqnarray}
where $E[\hat{\theta}^{add,RD}(\textbf{x}^*)]=\theta^{add}(\textbf{x}^*)+O(h^{s-1})$ and the standardizing factor
\begin{eqnarray*}
V_3 &=&  \frac{1}{n(2\pi)^{2d}}\int_{\mathbb{R}^d} \bigg(  \int_{\mathbb{R}^d} e^{i\langle \textbf{w},\textbf{s}\rangle } \Big( \sum_{j=1}^m e^{-i\langle \textbf{w}_{I_j},\textbf{x}_{I_j}^*\rangle }L_{I_j^c}(\textbf{w}_{I_j^c} )\Big) \frac{\Phi_K(h\textbf{w})}{\Phi_\psi(\textbf{w})} d\textbf{w} \bigg)^2 \frac{(\sigma^2+g(\textbf{s})^2)f(\textbf{s})}{\max\{ f(\textbf{s}),f(\frac{1}{\textbf{a}_n})\}}d\textbf{s}.
\end{eqnarray*}
satisfies
\begin{eqnarray*}
C_ln^{1/2}h^{d/2+\beta-\gamma_{min}}f(\textbf{a}_\textbf{n}^{-1})^{1/2}\le V_3^{-1/2} \le C_u n^{1/2}h^{d/2+\beta-\gamma_{min}}.
\end{eqnarray*}
\end{theorem}

\subsection{Additive estimator for fixed design}

The asymptotic properties of the additive estimator $\hat{\theta}^{add,RD}$ defined in \eqref{margint} under the fixed design assumption have been studied by \cite{birbishil2012} and in this Section we investigate the asymptotic properties of the alternative estimator defined in Section \ref{sec22}. Our first result, Lemma \ref{lemfd}, gives the weak convergence of $\hat{\theta}_{I_j}^{FD}$, whereas Theorem \ref{theofd} contains the asymptotic distribution of the estimator $\hat{\theta}^{add,FD}$ defined in \eqref{est}. The proofs  are again deferred to Section \ref{sec5}.

\begin{lemma}
\label{lemfd}
Consider the inverse regression model under the fixed design assumption (FD). Let Assumptions {\rm \ref{asspsiadd}(B), \ref{assk}, \ref{asstheta}(C) \textit{and} \ref{asstheta}(D)}  be fulfilled for some $j \in \{1,...,m\}$, $h \rightarrow 0$ and $a_n \rightarrow 0$ as $n \rightarrow \infty$ such that
 \begin{eqnarray*}
 n^dh^{d_j+2\beta_j}a_n^{d_j} \rightarrow \infty \quad \mbox{and} \quad  n^2h^{2+d_j+\beta_j}a_n^3 \rightarrow \infty,
\end{eqnarray*}
 then
\begin{eqnarray}
U_{nj}(\textbf{x}_{I_j}^*)^{-1/2} ( \hat{\theta}_{I_j}^{FD}(\textbf{x}_{I_j}^*)-E[\hat{\theta}_{I_j}^{FD}(\textbf{x}_{I_j}^*)] )  \Rightarrow {\cal{N}}(0,1), \label{weaklemFD}
\end{eqnarray}
where the normalizing sequence is defined by
\begin{eqnarray*}
U_{nj}(\textbf{x}_{I_j}^*)&=& \frac{\sigma^2}{(2n+1)^{d-d_j}} \sum_{k_{I_j} \in \{-n,...,n\}^{d_j}} w_{k_{I_j},n}(\textbf{x}_{I_j}^*)^2,
\end{eqnarray*}
the weights $w_{k_{I_j},n}$ are defined in \eqref{wfdadd} and
\[E[\hat{\theta}_{I_j}^{FD}(\textbf{x}_{I_j}^*)] = \theta_{I_j}^{add}(\textbf{x}_{I_j}^*)+O(h^{s-1})+O(n^{-2}h^{-d_j-\beta_j-2}a_n^{-3}).\]
\end{lemma}
The result of Theorem \ref{theofd} below follows immediately from Lemma \ref{lemfd}. The bias is of the same order as the bias in Lemma \ref{lemfd} and we define $j^*= {\rm argmax}_j \: (d_j+2\beta_j)$.

\begin{theorem}
\label{theofd}
Consider the inverse regression model under the fixed design assumption (FD). Let Assumptions {\rm \ref{asspsiadd},\ref{assk}, \ref{asstheta}(C) \textit{and} \ref{asstheta}(D)} be fulfilled, $h \rightarrow 0$ and $a_n \rightarrow 0$ as $n \rightarrow \infty$ such that
\[n^dh^{d_{j^*}+2\beta_{j^*}}a_n^{d_{j^*}} \rightarrow \infty \quad \mbox{and} \quad n^2h^{2+d_{j^*}+\beta_{j^*}}a_n^3 \rightarrow \infty. \]
 Then
\begin{eqnarray}
 U_n(\textbf{x}^*) ^{-1/2} ( \hat{\theta}^{add,FD}(\textbf{x}^*)-E[\hat{\theta}^{add,FD}(\textbf{x}^*)] ) \Rightarrow {\cal{N}}(0,1), \label{weak3}
\end{eqnarray}
where the normalizing sequence is defined by
\begin{eqnarray*}
U_n(\textbf{x}^*) &=& \sigma^2 \sum_{\textbf{k} \in \{-n,...,n\}^d} \left( \sum_{j=1}^m \frac{1}{(2n+1)^{d-d_j}} w_{k_{I_j},n}(\textbf{x}_{I_j}^*) \right)^2,
\end{eqnarray*}
the weights $w_{k_{I_j},n}$ are defined in \eqref{wfdadd} and
\[E[\hat{\theta}^{add,FD}(\textbf{x}^*)]=\theta^{add}(\textbf{x}^*)+O(h^{s-1})+O(\frac{1}{n^2h^{2+ \max_j d_j+ \max_j \beta_j}a_n^3}).\]
\end{theorem}

\begin{rem}
{\rm
\color{white}{a} \color{black}
\begin{enumerate}
\item[(1)] The normalizing sequence $U_n(\textbf{x}^*)$ in \eqref{weak3} is of order $n^dh^{d_{j^*}+2\beta_{j^*}}a_n^{d_{j^*}}$.
\item[(2)] The bias of the additive estimator in the fixed design case is only vanishing if the subsets $I_j$ in the decomposition \eqref{thetaadd} satisfy $d_j\le 3$ for all $j=1,...,m$.
\item[(3)] Theorem 3.2 can easily be extended to multiplicative models of the form \eqref{model} with
\begin{eqnarray*}
\theta(\textbf{x}^*) &=& \prod_{j=1}^m \theta_{I_j}(\textbf{x}_{I_j}^*)
\end{eqnarray*}
if the convolution function $\psi$ is also multiplicative. Otherwise the estimator is not consistent and other techniques such as the marginal integration method have to be used.
\end{enumerate}}
\end{rem}

\section{Dependent data} \label{sec3a}
\def\theequation{4.\arabic{equation}}
\setcounter{equation}{0}
In this Section we briefly discuss the case of dependent data. To be precise
we assume that the errors in the inverse regression models have an MA($q$) structure. Under the random design
assumption this structure is given by
\begin{equation} \label{marandiom}
\epsilon_t = Z_t+\beta_1Z_{t-1}+...+\beta_qZ_{t-q},
\end{equation}
 where $\{Z_t,\}_{t \in \mathbb{Z}}$ denotes a white noise process with variance $\sigma^2$. A careful inspection of the proof of Theorem \ref{theoaddrd}, which is based
 on the investigation of the asymptotic properties of cumulants shows that the result of Theorem \ref{theoaddrd} remains valid
 under this assumption.

\begin{theorem} ~\\
\label{theomaqrd}
(1) Consider the inverse regression model \eqref{model2} under the random design assumption (RD). If
the Assumptions of Theorem \ref{theoaddrd} are satisfied, then
\begin{eqnarray}
 V_1 ^{-1/2} \big(\hat{\theta}^{RD}(\textbf{x}^*)-E[\hat{\theta}^{RD}(\textbf{x}^*)]\big) \Rightarrow {\cal{N}}(0,1), \label{weak1}
\end{eqnarray}
where the normalizing sequence is given by
\begin{eqnarray*}
V_1 &=& \frac{1}{nh^{d}(2\pi)^{2d}}\int_{\mathbb{R}^d}\Big( \int_{\mathbb{R}^d} e^{-i\langle \textbf{s},(\textbf{x}^*/h-\textbf{y})\rangle }\frac{\Phi_K(\textbf{s}) }{\Phi_{\psi}(\frac{\textbf{s}}{h})} d\textbf{s} \Big)^2\frac{(\sigma^2\sum_{k,l=0}^q  \beta_k\beta_{l}+g^2(h\textbf{y}))f(h\textbf{y})}{\max\{f(h\textbf{\textbf{y}}),f(\frac{1}{\textbf{a}_n})\}^2}d\textbf{y},
\end{eqnarray*}
$\beta_0=1$ and $E[\hat{\theta}^{RD}(\textbf{x}^*)]=\theta(\textbf{x}^*)+O(h^{s-1})$. \\
(2)  If the assumptions of Theorem \ref{theoremaddrd} are satisfied, then the appropriately standardized additive estimator
$\hat{\theta}^{add,RD}$ converges weakly to a standard normal distribution, that is
\begin{eqnarray}
V_3^{-1/2} \big( \hat{\theta}^{add,RD}(\textbf{x}^*)-E[\hat{\theta}^{add,RD}(\textbf{x}^*)] \big) \Rightarrow {\cal{N}}(0,1), \label{weak2}
\end{eqnarray}
where the standardizing factor is given by
\begin{eqnarray*}
V_3 &=&  \frac{1}{n(2\pi)^{2d}}\int_{\mathbb{R}^d} \bigg(  \int_{\mathbb{R}^d} e^{i\langle \textbf{w},\textbf{s}\rangle } \Big( \sum_{j=1}^m e^{-i\langle \textbf{w}_{I_j},\textbf{x}_{I_j}^*\rangle }L_{I_j^c}(\textbf{w}_{I_j^c} )\Big) \frac{\Phi_K(h\textbf{w})}{\Phi_\psi(\textbf{w})} d\textbf{w} \bigg)^2 \frac{(\sigma^2\sum_{k,l=0}^q  \beta_k\beta_{l}+g(\textbf{s})^2)f(\textbf{s})}{\max\{ f(\textbf{s}),f(\frac{1}{\textbf{a}_n})\}}d\textbf{s}.
\end{eqnarray*}
and $E[\hat{\theta}^{add,RD}(\textbf{x}^*)]=\theta^{add}(\textbf{x}^*)+O(h^{s-1})$.
\end{theorem}

Under the assumption of a fixed design on a grid we consider an error process with an MA($q$) structure defined by
\begin{eqnarray}
\epsilon_\textbf{k} &=& \sum_{\textbf{r} \in \{-q,...,q\}^d} \beta_\textbf{r} Z_{\textbf{k}-\textbf{r}}, \label{MAq}
\end{eqnarray}
where $\{Z_\textbf{j}\}_{j\in \mathbb{Z}^d}$ are i.i.d. random variables with mean zero and variance $\sigma^2$.
This means, that the noise terms are influenced by all shocks, which have a distance on the lattice lower or equal
$q$ regarding the $\infty$-norm.  The following result can be obtained by similar arguments as used for the proof of Theorem \ref{theofd}.

\begin{theorem}
\label{lemfd2}
Consider the inverse regression model \eqref{model3} under the fixed design assumption with an MA(q) dependent error process. If the assumptions of Lemma \ref{lemfd} are satisfied we have
\begin{eqnarray*}
V_{MA}^{-1/2}(\textbf{x}^*) \left( \hat{\theta}^{add,FD}(\textbf{x}^*)-E[\hat{\theta}^{add,FD}(\textbf{x}^*)]\right) \Rightarrow \mathcal{N}(0,1)
\end{eqnarray*}
where the normalizing sequence is given by
\begin{eqnarray*}
V_{MA}(\textbf{x}^*)&=&\sigma^2\sum_{ \substack{\textbf{l} \in  \mathbb{Z}^d \\ \parallel \textbf{l} \parallel_\infty \le 2q}}\sum_{\textbf{r}_1 \in \{-q,...,q\}^d}\beta_{\textbf{r}_1}\beta_{\textbf{l}+\textbf{r}_1} \sum_{\textbf{k} \in \{-n,...,n\}^{d}} |\sum_{j=1}^m \frac{1}{(2n+1)^{d-d_j}}w_{\textbf{k}_{I_j},n}(\textbf{x}_{I_j}^*)|^2.
\end{eqnarray*}
and $E[\hat{\theta}^{add,FD}(\textbf{x}^*)]=\theta^{add}(\textbf{x}^*)+O(h^{s-1})+O(\frac{1}{n^2h^{2+ \max_j d_j+ \max_j \beta_j}a_n^3})$.
\end{theorem}

\begin{rem} \label{MAunendlich} {\rm
If $\epsilon_t$ has an MA($\infty$) representation Theorem \ref{theomaqrd} and \ref{lemfd2} will not hold in general, because without additional assumptions the $l$-th cumulant of the normalized statistic does not converge to zero for all $l \ge 3$. }
\end{rem}

\section{Finite sample properties}  \label{sec4}
\begin{figure}[hbt]
\center
  \includegraphics[width = 8cm,height=6cm]{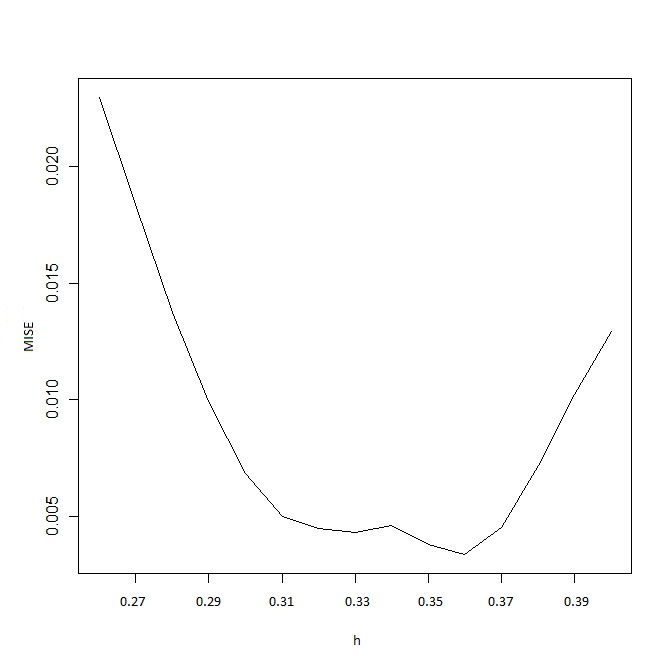}
\caption{\textit{\label{fig11} MISE of the estimator  $\hat{\theta}^{add,FD}$ for different bandwidths in model \eqref{sig1}, where $\sigma = 0.5$}  }
\hfill
\end{figure}

\begin{figure}[htb!]
	\includegraphics[width=8.cm]{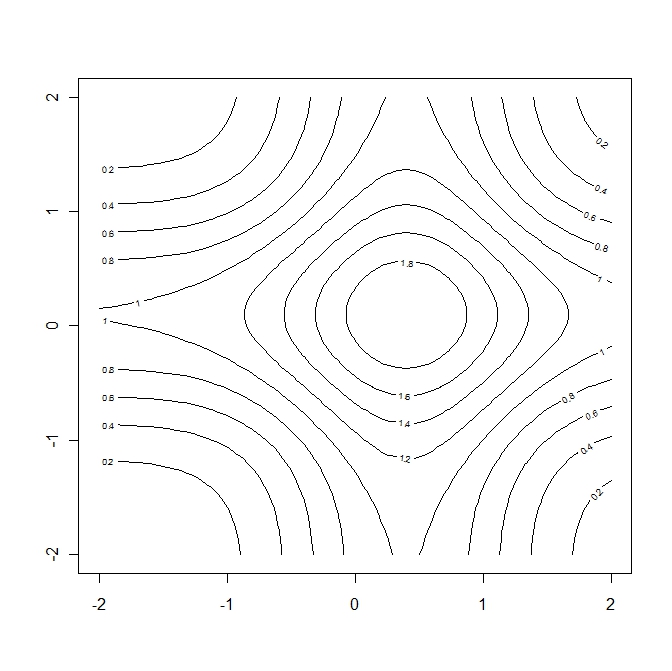}
	\includegraphics[width=8.cm]{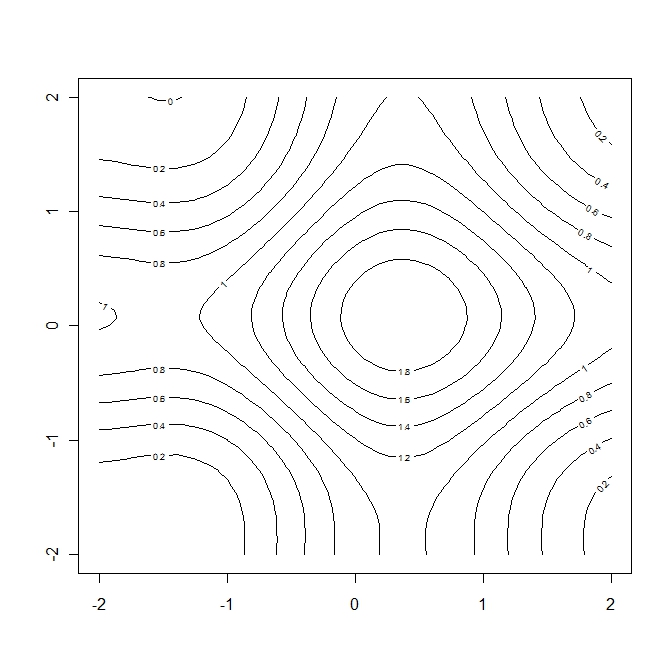}
	\includegraphics[width=8.cm]{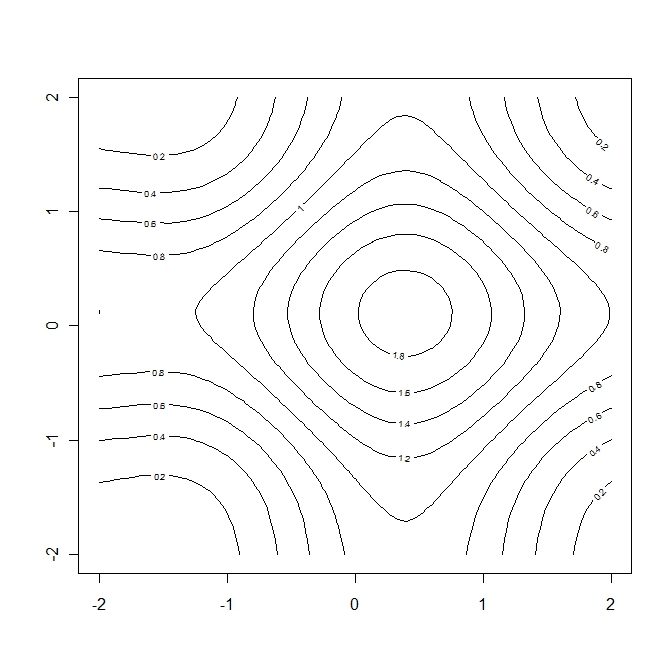} \hfill
	\includegraphics[width=8.cm]{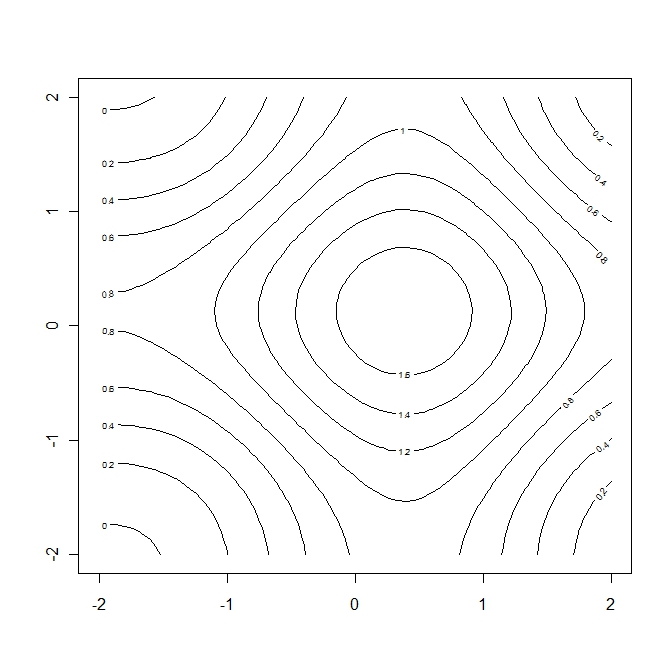}
\caption{{\it
\label{fig2}
Contour plot of the function $\theta^{(1)}$ defined in \eqref{sig1} (left upper panel) and its estimates $\hat{\theta}^{add,FD}$ defined in \eqref{est} with different bandwidths. Upper right panel:  $h=0.32$; Lower left panel: $h=0.36$ (which minimizes the MISE); Lower right panel: $h=0.4$; } }
\end{figure}

\def\theequation{5.\arabic{equation}}
\setcounter{equation}{0}
In this Section we investigate the finite sample properties of the new estimators and also provide a comparison with competing methods. We first investigate the case of a fixed design in model \eqref{model} with the convolution function
\begin{eqnarray*}
\psi(x_1,x_2) &=& \frac{9}{4}e^{-3(|x_1|+|x_2|)},
\end{eqnarray*}
 and two additive signals
\begin{eqnarray}
\theta^{(1)}(x_1,x_2) &=& e^{-(x_1-0.1)^2}+e^{-(x_2-0.4)^2} \label{sig1}
\\ \theta^{(2)}(x_1,x_2) &=& e^{-|x_1-0.4|}+2e^{-2x_2^2} \label{sig2}
\end{eqnarray}
For the kernel $K$ in the Fourier transform $\Phi_K$ we use the kernel $K(\textbf{x})=\frac{\sin(x_1)\sin(x_2)}{\pi^2x_1x_2}$. We consider a fixed design on the grid $\{(\frac{k_1}{na_n},\frac{k_2}{na_n} \: | \: k_1,k_2 \in \{-n,...,n\}\}$ with $N=(2n+1)^2$ points where $n \in \{30,50\}$. In both cases we choose the design parameter as $a_n=0.25$, such that the cube $[\frac{-1}{a_n},\frac{1}{a_n}]^2$ covers most of the region where the functions $\theta^{(1)}$ and $\theta^{(2)}$ deviate significantly from 0. In all simulations we use (independent) noise terms, which are normal distributed with mean 0 and variance 0.25.
\\ The bandwidth $h$ in the estimator \eqref{estcomp} is chosen such that the mean integrated squared error (MISE)
\begin{eqnarray*}
\mathbb{E} \Bigr[\int_{\mathbb{R}^2} (\hat{\theta}(\textbf{x})-\theta(\textbf{x}))^2 d\textbf{x}\Bigr]
\end{eqnarray*}
is minimized. Figure \ref{fig11} shows a typical example of the MISE as a function of the bandwidth $h$.
Figure \ref{fig2}  shows the contour plot of the function $\theta^{(1)}$ defined in \eqref{sig1} and contour plots of three typical additive estimates where $n=50$ and the bandwidths are chosen as $h=0.32,0.36,0.4$  (the bandwidth $h=0.36$ minimizes the MISE). We observe that the shapes in all figures are very similar.  The bandwidths $h=0.32$ and $h=0.4$ yield stronger deviations from the true function especially at the boundary, but the main structure is even for these choices still recovered. Because other simulations showed a similar picture we conclude that small changes in the bandwidth do  not effect the general structure of the estimator significantly.
\\In order to investigate the finite sample properties of the
new estimate $\hat{\theta}^{add,FD}$ defined in \eqref{est}  we performed 1000 iterations with the signal $\theta^{(2)}$ (the results for the signal $\theta^{(1)}$ are similar and are not depicted for the sake of brevity). The simulated mean, variance and mean squared error (MSE) of $\hat{\theta}^{add,FD}$ are given in Table \ref{tab1}
  for different choices of $\textbf{x}=(x_1,x_2)$ where the sample size is $N=10201$ and the variance of the errors is $0.25$. We observe that 
  in most cases the mean squared error is dominated by
  the bias.
\medskip

\begin{table}[!htb]
\centering
{\scriptsize
\begin{tabular}{|| c|c|c|c|c|c|c ||}
\hline
 $N$  &$x_1$& $x_2$& $\theta^{(2)}(\textbf{x})$ &  $E[\hat{\theta}(\textbf{x})]$ & Var($\hat{\theta}(\textbf{x})$)&MSE($\hat{\theta}(\textbf{x}$) \\ \hline
  &    &  -1.6&  0.1473  & 0.2522 &0.0017  & 0.0127   \\
  &     & -0.8 &  0.3131 & 0.3805  & 0.0017 & 0.0063 \\
 10201 & -1.6  & 0 & 0.6823  & 0.8296  &0.0017  & 0.0234  \\
  &     & 0.8 & 0.6823    &  0.8159 & 0.0017 & 0.0195 \\
  &     & 1.6 & 0.3131  & 0.3827  &0.0017  & 0.0065\\ \hline
 &     & -1.6 &  0.6914 & 0.8216  &0.0017  &0.0187 \\
  &    &-0.8  &  0.8573  & 0.9446  & 0.0018 &0.0094  \\
 10201  &  -0.8 & 0 & 1.2264    &1.3977 & 0.0017 & 0.0310  \\
  &     &0.8  & 1.2264  & 1.3864  &0.0017  & 0.0273 \\
 &     & 1.6 &  0.8573  & 0.9496 &0.0018  & 0.0103  \\ \hline
 &     &-1.6  & 2.1353  & 2.1887  &0.0018  & 0.0046 \\
  &     & -0.8 &  2.3012  &  2.3123 &0.0017 &0.0018  \\
 10201  &  0 & 0 & 2.6703   & 2.7640&0.0018  & 0.0106 \\
  &     & 0.8 & 2.6703  &  2.7548 &0.0016  &0.0087 \\
 &     & 1.6 &2.3012    &2.3178  & 0.0018 & 0.0020  \\ \hline
 &     &-1.6  &  0.6914 & 0.8181  & 0.0017 &0.0178 \\
  &    &-0.8  & 0.8573   & 0.9445  & 0.0018 & 0.0094 \\
 10201  &  0.8 & 0 & 1.2264   & 1.3967& 0.0017 & 0.0307 \\
  &     &0.8  &  1.2264 & 1.3864 & 0.0017 &0.0273 \\
 &     &1.6  & 0.8573   &0.9496 & 0.0018 & 0.0103  \\ \hline
 &     & -1.6 &  0.1473 & 0.2532  & 0.0016 &0.0128 \\
  &     & -0.8 &   0.3131 &  0.3785 & 0.0017 & 0.0060  \\
 10201  &  1.6 & 0 & 0.6823   &0.8290 &0.0018  &0.0233  \\
  &     & 0.8 &0.6823   &0.8168  &  0.0019& 0.0200\\
 &     & 1.6 &  0.3131  & 0.3855 & 0.0017 & 0.0069   \\ \hline
\end{tabular}
}
\caption{{\it \label{tab1}
Mean, variance and mean squared error of the new additive estimator $\hat{\theta}=\hat{\theta}^{add,FD}$ in the case of a fixed design. The model is given by \eqref{sig2} with variance $\sigma^2=0.25$.}}
\end{table}
In the second part of this section we compare three different estimates for the signal in the inverse regression model \eqref{model}. The first estimate for $\theta$ is the statistic $\hat{\theta}^{add,FD}$ proposed in this paper [see formula \eqref{est}]. The second method is the marginal integration estimator suggested by \cite{birbishil2012} and the third method is the non additive estimate of \cite{bb2008}. The results are shown in Table  \ref{tab2}
for the sample size $N=3721$ and selected values of the predictor. We observe that the additive estimate of \cite{birbishil2012} improves the unrestricted estimate with respect to mean squared error by 20-50\%. However, the new additive estimate $\hat{\theta}^{add,FD}$ yields a much larger improvement. The MSE is about 14 and 7-10 times smaller than the MSE obtained by the unrestricted estimator  or the estimator proposed by \cite{birbishil2012}. Further simulations for the signal $\theta^{(2)}$ in \eqref{sig2} show similar results and not depicted for the sake of brevity.

\begin{table}[!htb]
\centering
{\scriptsize
\begin{tabular}{|| c|c|c|c|c|c|c|c ||}
\hline
& $N$  &$x_1$& $x_2$& $\theta^{(1)}(\textbf{x})$ &  $\mathbb{E}\hat{\theta}(\textbf{x})$ & Var $\hat{\theta}(\textbf{x})$&MSE $\hat{\theta}(\textbf{x})$ \\ \hline
  & 3721&0  &0    &  1.8422  & 1.9667 & 0.0516 & 0.0671  \\
$\hat{\theta}^{RD}$ & 3721 & 0 & 1    &   1.6877  & 1.6983  & 0.0458 & 0.0459 \\
& 3721 & 1 &  1   &   1.1425 &  1.1909 & 0.0329 & 0.0352 \\
&3721 & 1 & 1.8    & 0.5857   & 0.6624& 0.0189 & 0.0248  \\ \hline
  & 3721&0  &0    &  1.8422  & 1.8680 &0.0440  & 0.0301   \\
$\hat{\theta}^{add,RD}$ & 3721 & 0 & 1    &   1.6877  & 1.6405   & 0.0195 & 0.0217  \\
& 3721 & 1 &  1   &   1.1425 & 1.3371  & 0.0232 &  0.0610 \\
&3721 & 1 & 1.8    & 0.5857   & 0.8184 & 0.0199 &  0.0740 \\ \hline \hline
 &3721  & 0 &0 &  1.8422 & 1.8123 & 0.0426 & 0.0435 \\
$\hat{\theta}^{FD}$&3721  &  0 &  1&  1.6877 & 1.7305  & 0.0425 &0.0443 \\
& 3721  & 1  &1  &  1.1425 &  1.2143 & 0.0418 & 0.0470\\
&3721  &  1 & 1.8 &  0.5857 &  0.4774 & 0.0416 & 0.0533  \\ \hline
  & 3721&0  &0    &  1.8422  & 1.8234 & 0.0027 & 0.0031  \\
$\hat{\theta}^{add,FD}$ & 3721 & 0 & 1    &   1.6877  & 1.6589  & 0.0024 & 0.0032 \\
& 3721 & 1 &  1   &   1.1425 &  1.1097 & 0.0025 & 0.0036 \\
&3721 & 1 & 1.8    & 0.5857   & 0.5494 & 0.0023 & 0.0036  \\ \hline
&  3721  &  0 & 0 & 1.8422  & 1.8874 & 0.0194 & 0.0214 \\
$\hat{\theta}^{BBH}$& 3721  &  0 & 1 & 1.6877 & 1.7316 & 0.0191 & 0.0210\\
& 3721  &  1 & 1 &  1.1425  &1.1833 & 0.0201 & 0.0218 \\
&3721  &  1 & 1.8 &  0.5857 &0.4438 & 0.0207 & 0.0408  \\  \hline

\end{tabular}
}
\caption{{\it \label{tab2}
Mean, variance and mean squared error of the unrestricted estimator $\hat{\theta}^{FD}$ proposed in \cite{bb2008}, the estimator $\hat{\theta}^{BBH}$ proposed by \cite{birbishil2012}  and the new  estimators $\hat{\theta}^{RD}$, $\hat{\theta}^{add,RD}$  and  $\hat{\theta}^{add,FD}$ proposed in this paper. The model is given by \eqref{sig1}, where $\sigma^2=0.25$.}}
\end{table}

For the sake of comparison, the first two rows of Table \ref{tab2} contain results of the estimators $\hat \theta^{RD}$ and $\hat \theta^{add,RD}$, where the explanatory variables follow a uniform distribution on the same cube $[\frac{1}{a_n},\frac{1}{a_n}]^2$ as used for the fixed design. We observe a similar behaviour of the unrestricted estimators  under the fixed and random design assumption. This corresponds to the asymptotic theory, which shows that in the case of a uniform distribution the unrestricted estimators converge with the same rate of convergence (see Remark \ref{rem3.2}). On the other hand, the additive estimator $\hat \theta^{add,RD}$ produces a substantially larger mean squared error compared to the estimator $\hat \theta^{add,FD}$, which is of similar size as the mean squared error of the estimator proposed by \cite{birbishil2012}.

\begin{figure}[hbt!]
	\includegraphics[width=8cm]{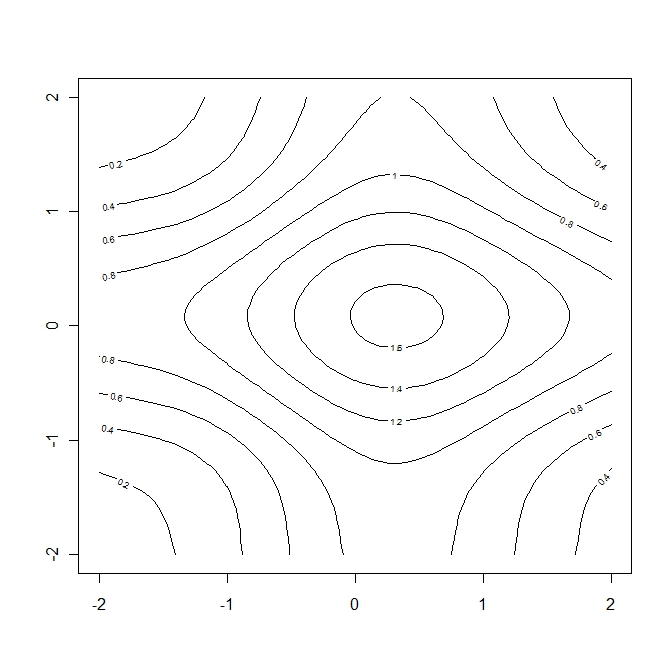}
	\includegraphics[width=8cm]{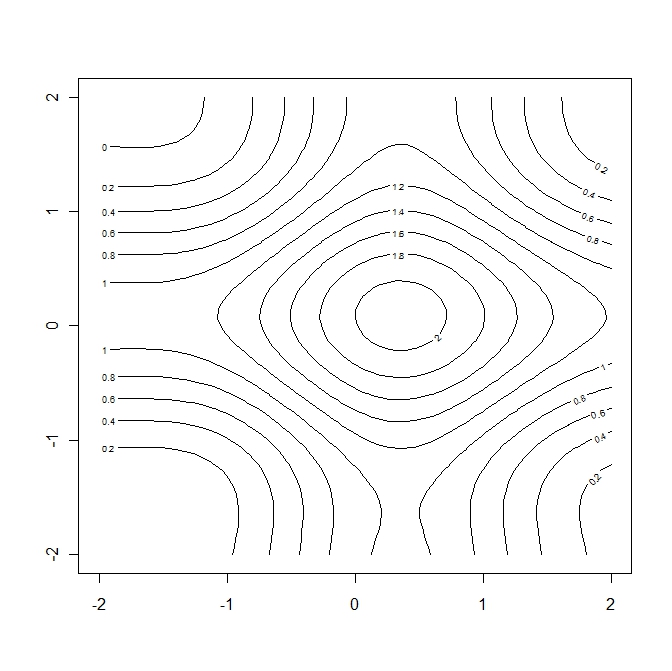}
	\includegraphics[width=8cm]{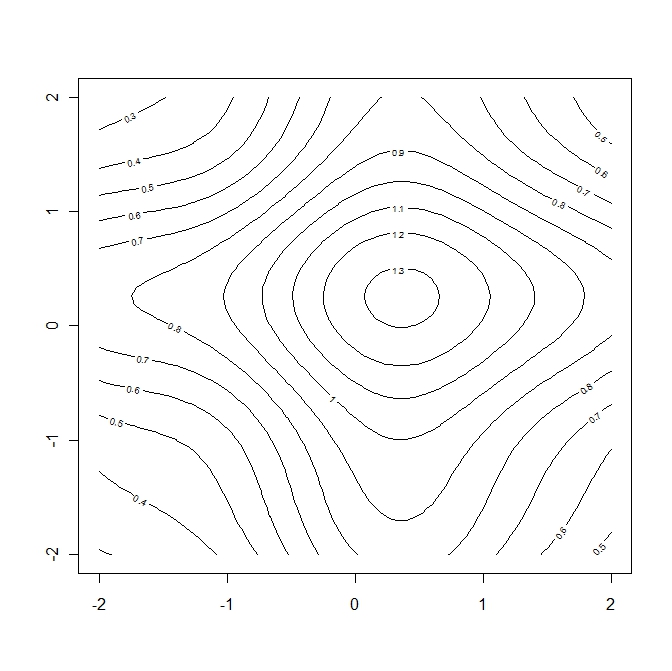} \hfill
	\includegraphics[width=8cm]{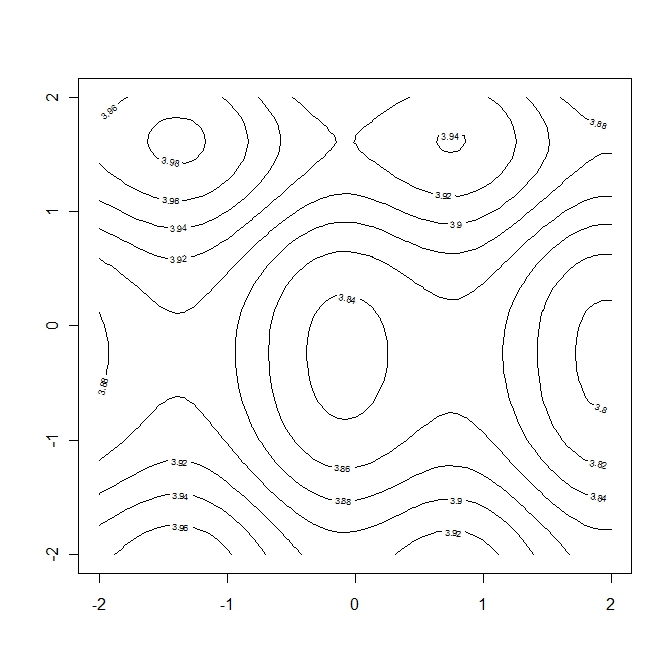}	
\caption{{ \it
\label{fig3}
Contour plot of the estimate $\hat{\theta}^{add,FD}$ of $\theta^{(1)}$ with misspecificated convolution function. Upper left panel: $\psi$ misspecificated as  Lap(0,$\frac{1}{3})$, where the true convolution function is Lap(0,1); Upper right panel: $\psi$ misspecificated as  Lap(0,$\frac{1}{3})$, where the true convolution function is Lap(0,$\frac{1}{5}$) ; Lower left panel: $\psi$ misspecificated as  Lap(0,$\frac{1}{3})$, where the true convolution function is ${\cal N}(0,1)$; Lower right panel: $\psi$ misspecificated as Lap(0,$\frac{1}{3})$,  where the true convolution function is  Gamma(3,2). The model is given by \eqref{sig1}, where $\sigma^2=0.25$.}}
\end{figure}

Because the performance of the estimators depends on the correct specification of the convolution function $\psi$ we next investigate the performance of the estimators under misspecification of the function $\psi$. In Figure \ref{fig3} we display the contour plots of the estimates $\hat{\theta}^{add,FD}$, where in every panel the convolution function is misspecificated as Laplace distribution $Lap(\alpha,\beta)$ with parameters $\alpha=0$ and $\beta=\frac{1}{3}$. In the upper left and upper right panel the $\beta$ parameter of the Laplace distribution $Lap(\alpha,\beta)$ is misspecificated, whereas in the lower left panel the true convolution function is the density of a standard normal distribution and in the lower right panel it is a gamma distribution. We observe, that a miss-specification of the shape of the convolution function (as it occurs if a Laplace density is used instead of the density of a Gamma(3,2) distribution) yields to an estimator with a different structure as the true signal (see the lower right panel in Figure 3).  All other panels show the same structure as the upper left panel Figure \ref{fig2} which gives the contour plot of the true signal $\theta^{(1)}$. This indicates that the structure of the signal can be reconstructed, as long as the chosen convolution kernel exhibits similar modal properties as the  ``true kernel'' . However, we also observe from
Figure \ref{fig3}
 that the levels of the contour differ from those of  the true signal.

We conclude this section with a brief discussion of the performance of the unrestricted estimator $\hat \theta^{RD}$ under the assumption (RD) of a non-uniform random design. In Table \ref{tab3} we display
the simulated mean, variance and mean squared error for various distributions of the predictor $\bf{X}$, where the components are independent  and identically distributed.
In most cases we observe  similar results for the bias, independently of the distribution of $\bf{X}$ and the choice of the sequence $a_n$. On the other hand the mean squared error is dominated
by the variance, which depends sensitively on the choice of the  parameter  $a_n$. This observation corresponds  with the representation of the asymptotic variance of  $\hat \theta^{RD}$
   in formula \eqref{v1} of Theorem \ref{theoaddrd}. We also observe that the impact of the
distribution of the explanatory variable on the variance of the estimate  $\hat \theta^{RD}$   is much smaller.
\begin{table}[!htb]
\centering
{\scriptsize
\begin{tabular}{|| c|c|c|c|c|c|c|c|c ||}
\hline
${\bf X}$& $N$ &$a_n$ &$x_1$& $x_2$& $\theta^{(1)}(\textbf{x})$ &  $\mathbb{E}\hat{\theta}(\textbf{x})$ & Var $\hat{\theta}(\textbf{x})$&MSE $\hat{\theta}(\textbf{x})$ \\ \hline
  & 10201&0.25  &0 &0   &  1.8422  & 1.7421 & 0.0297 & 0.0397 \\
$U[\frac{-1}{a_n},\frac{1}{a_n}]$ & 10201&0.25 & 0 & 1    &   1.6877  & 1.7163  & 0.0272 & 0.0283 \\
& 10201&0.25 & 1 &  1   &   1.1425 &  1.2858 & 0.0194 & 0.0399 \\
&10201 &0.25& 1 & 1.8    & 0.5857   & 0.6105& 0.0117& 0.0123  \\ \hline
  & 10201&0.5&0  &0    &  1.8422  & 1.4957 &0.0076  & 0.1277   \\
$U[\frac{-1}{a_n},\frac{1}{a_n}]$ & 10201&0.5 & 0 & 1    &   1.6877  & 1.8123   & 0.0070 & 0.0225  \\
& 10201&0.5 & 1 &  1   &   1.1425 & 1.5438  & 0.0044 &  0.1654 \\
&10201 &0.5& 1 & 1.8    & 0.5857   & 0.5695 & 0.0023 &  0.0026 \\ \hline
  & 10201&0.25&0  &0    &  1.8422  & 1.8512  & 0.3271  & 0.3271  \\
$N(0,1)$ & 10201&0.25 & 0 & 1    &   1.6877  & 1.7019  & 0.7098 & 0.7100 \\
& 10201&0.25 & 1 &  1   &   1.1425 &  1.2038 & 0.7077 & 0.7115 \\
&10201 &0.25& 1 & 1.8    & 0.5857   & 0.5983 & 0.4477 & 0.4479  \\ \hline
&  10201 &0.5 &  0 & 0 & 1.8422  & 1.8229 & 0.0079 & 0.0083\\
$N(0,1)$& 10201 &0.5 &  0 & 1 & 1.6877 & 1.7466 & 0.0107 & 0.0143\\
& 10201 &0.5 &  1 & 1 &  1.1425  &1.2531 & 0.0114 & 0.0236 \\
& 10201 &0.5 &  1 & 1.8 &  0.5857  &0.6366 & 0.0135 & 0.0161 \\ \hline
  & 10201&0.25&0  &0    &  1.8422  & 1.8758 &0.0174  & 0.0185   \\
$t(2)$ & 10201&0.25 & 0 & 1    &   1.6877  & 1.7129   & 0.0255 & 0.0261 \\
& 10201&0.25& 1 &  1   &   1.1425 & 1.1786  & 0.0271 &  0.0284\\
&10201 &0.25& 1 & 1.8    & 0.5857   & 0.6138 & 0.0324 &  0.0332 \\ \hline
 &10201 &0.5 & 0 &0 &  1.8422 & 1.8590 & 0.0115 & 0.0118 \\
$t(2)$&10201 &0.5 &  0 &  1&  1.6877 & 1.7260  & 0.0158 &0.0173 \\
& 10201 &0.5 & 1  &1  &  1.1425 &  1.2069 & 0.0182 & 0.0223\\
&10201 &0.5 &  1 & 1.8 &  0.5857 &  0.6275 & 0.0174 & 0.0191  \\ \hline
\end{tabular}
}
\caption{{\it \label{tab3} Mean, variance and mean squared error of the unrestricted estimator $\hat{\theta}^{RD}$ proposed in this paper
for different distributions of the explanatory  variables ${\bf X}$ and different choices   for the parameter $a_n$. The model is given by  \eqref{sig1} and the variance is
$\sigma^2=0.25$. }}
\end{table}

\medskip

{\bf Acknowledgements.}
This work has been supported in part by the
Collaborative Research Center ``Statistical modeling of nonlinear
dynamic processes'' (SFB 823, Teilprojekt C1, C4) of the German Research Foundation
(DFG).

\bibliographystyle{apalike}

\bibliography{hilbisdet}

\section{Appendix}  \label{sec5}
\def\theequation{5.\arabic{equation}}
\setcounter{equation}{0}

For the proofs we make frequent use of the cumulant method,  which is a common tool in time series analysis. Following \cite{bril2001} the $r$-th order joint cumulant $cum(Y_1,...,Y_r)$ of a $r$-dimensional complex valued random vector $(Y_1,...,Y_r)$ is given by
\begin{eqnarray}
cum(Y_1,...,Y_r) &=&  \sum (-1)^{p-1} (p-1)!\Big(E\prod_{j \in \nu_1} Y_j\Big) \dots \Big(E\prod_{j \in \nu_p} Y_j\Big), \label{cum}
\end{eqnarray}
where  we assume the existence of moments of order $r$, i.e. $E( |Y_j^r|) < \infty$  $(j=1,...,r)$ and  the summation extends over all partitions $(\nu_1,...,\nu_p),p=1,...,r$ of $(1,...,r)$.
If we choose $Y_j=Y,j=1,...,r$ we denote with
$cum_r(Y)=cum(Y,...,Y)$
 the $r$-th order cumulant of a univariate random variable. The following  properties of the cumulant will be used frequently in our proofs [see e.g. \cite{bril2001}].
\begin{enumerate}
\item[(B1)] $cum(a_1Y_1,...,a_rY_r) = a_1\dots a_rcum(Y_1,...,Y_r)$ \quad for constants  $a_1,...,a_r\in \mathbb{C}$
\item[(B2)] if any group of the Y's is independent of the remaining Y's, then  $cum(Y_1,...,Y_r)=0$
\item[(B3)] for the random variable $(Z_1,Y_1,...,Y_r)$ we have \[cum(Z_1+Y_1,Y_2,...,Y_r)=cum(Z_1,Y_2,...,Y_r)+cum(Y_1,Y_2,...,Y_r)\]
\item[(B4)] if the random variables $(Y_1,...,Y_r)$ and $(Z_1,...,Z_r)$ are independent, then \[ cum(Y_1+Z_1,...,Y_r+Z_r) = cum(Y_1,...,Y_r) + cum(Z_1,...,Z_r) \]
\item[(B5)] $cum(Y_j) = E(Y_j)$ for $j=1,...,r$
\item[(B6)] $cum(Y_j,\overline{Y}_j) = Var(Y_j)$  for $j=1,...,r$
\end{enumerate}
We finally state a result which can easily be proven by using the definition \eqref{cum} and the properties of the mean.
\begin{theorem}
Let $\textbf{Y}=(Y_1,...,Y_r)$ be a random variable, $b_n$ a sequence and $C >0$ a constant with
\begin{eqnarray*}
E\Big[\prod_{j=1}^l|Y_{i_j}|\Big] &\le& C^lb_n^l \qquad \mbox{for all } 1\le l\le r,
\end{eqnarray*}
then
$
|cum(Y_{i_1},...,Y_{i_m})| \le (m-1)!C^mb_n^m \sum_{j=1}^m S_{m,j},
$
where $S_{m,j}$ denotes the Sterling number of the second kind.
\end{theorem}
 We will also make use of the fact that the normal distribution with mean $\mu$ and variance $\sigma^2$ is characterized by its cumulants, where the first two cumulants are equal to $\mu$ and $\sigma^2$ respectively and all cumulants of larger order are zero. To show asymptotic normality in our proofs we have to calculate the first two cumulants which give the asymptotic mean and variance and show in a second step that all cumulants of order $l\ge3$ are vanishing asymptotically. In the following discussion all constants which do not depend on the sample size (but may differ in different steps of the proofs) will be denoted by $C$.

\medskip

\textbf{Proof of Theorem \ref{theoaddrd}:}
For the sake of brevity we write $\hat{\theta}$ instead of $\hat{\theta}^{full,RD}$ throughout this proof. By the discussion of the previous paragraph we have to calculate the mean and the variance of $\hat{\theta}(\textbf{x}^*)$ and all cumulants of order $l\ge3$. We start with the mean conditional on $\textbf{X}=(\textbf{X}_1,...,\textbf{X}_n)$, which can be calculated as
\begin{eqnarray*}
E[\hat{\theta}(\textbf{x}^*)|\textbf{X}] &=&  \sum_{k=1}^n g(\textbf{X}_k)w_{n}(\textbf{x}^*,\textbf{X}_k)
\end{eqnarray*}
where the weights $w_{n}$ are defined in \eqref{wRD}. By iterative expectation we get
\begin{eqnarray*}
E[\hat{\theta}(\textbf{x}^*)] &=& \frac{1}{h^d(2\pi)^d} \int_{\mathbb{R}^d} g(\textbf{x}) \int_{\mathbb{R}^d} e^{-i\langle \textbf{s},(\textbf{x}^*-\textbf{x})\rangle /h}\frac{\Phi_K(\textbf{s}) }{\Phi_{\psi}(\frac{\textbf{s}}{h})}\frac{f(\textbf{x})}{\max\{f(\textbf{x}),f(\frac{1}{\bf{a}_n})\}} d\textbf{s}d\textbf{x},
\end{eqnarray*}
which yields a bias of the form $bias_{\hat{\theta}} = E[\hat{\theta}(\textbf{x}^*)]-\theta(\textbf{x}^*) = A_1+A_2,$ where (note that $\Phi_g=\Phi_\psi\cdot\Phi_\theta$)
\begin{eqnarray*}
A_1 &=& \frac{1}{h^d(2\pi)^d}\int_{\mathbb{R}^d} e^{-i\langle \textbf{s}, \textbf{x}^*\rangle /h}\Phi_K(\textbf{s})\Phi_\theta\left(\frac{\textbf{s}}{h}\right) d\textbf{s}-\theta(\textbf{x}^*)
\\ A_2 &=& \frac{1}{h^d(2\pi)^d}\int_{\mathbb{R}^d} e^{-i\langle \textbf{s}, \textbf{x}^*\rangle /h}\frac{\Phi_K(\textbf{s}) }{\Phi_{\psi}(\frac{\textbf{s}}{h})}\int_{\mathbb{R}^d} g(\textbf{x}) e^{i\langle s,\textbf{x}\rangle /h}\Big(\frac{f(\textbf{x})}{\max\{f(\textbf{x}),f(\frac{1}{\bf{a}_n})\}}-1\Big)d\textbf{x} d\textbf{s}
\end{eqnarray*}
 For the summand $A_1$ we can use exactly the same calculation as in \cite{bb2008} to obtain $A_1=O(h^{s-1})$. For the second term $A_2$ we have
\begin{eqnarray*}
A_2  &\le&\frac{1}{h^d(2\pi)^d}\int_{\mathbb{R}^d} \frac{|\Phi_K(\textbf{s})| }{|\Phi_{\psi}(\frac{\textbf{s}}{h})|}\int_{\mathbb{R}^d} |g(\textbf{x})| \Big|\frac{f(\textbf{x})}{\max\{f(\textbf{x}),f(\frac{1}{\bf{a}_n})\}}-1\Big|d\textbf{x} d\textbf{s}
\\ &\le& \frac{C}{h^{d+\beta}(2\pi)^d} \int_{([-\frac{1}{a_n},\frac{1}{a_n}]^d)^c} |g(\textbf{x})| \Big|\frac{f(\textbf{x})}{\max\{f(\textbf{x}),f(\frac{1}{\bf{a}_n})\}}-1\Big|d\textbf{x},
\end{eqnarray*}
where we used Assumption \ref{asspsiadd}(A) and 4 in the last inequality.
In the next step we will use the fact that $0\le\frac{f(\textbf{x})}{\max\left\{f(\textbf{x}),f(\frac{1}{\bf{a}_n})\right\}} \le1$  ($\textbf{x} \in \mathbb{R}^d$) and Assumption \ref{asstheta}(B) to obtain
\begin{eqnarray*}
A_2 &\le&  \frac{C}{h^{d+\beta}(2\pi)^d} \int_{([-\frac{1}{a_n},\frac{1}{a_n}]^d)^c} |g(\textbf{x})| \parallel \textbf{x} \parallel^r \frac{1}{\parallel \textbf{x} \parallel^r} d\textbf{x} = O\Big(\frac{a_n^r}{h^{d+\beta}}\Big) = O(h^{s-1}).
\end{eqnarray*}
This shows that the bias of $\hat{\theta}(\textbf{x}^*)$ is of order $ O(h^{s-1})$.
By the definition of $\hat{\theta}(\textbf{x}^*)$ and \eqref{wRD} it follows
\begin{eqnarray*}
V(\hat{\theta}(\textbf{x}^*)|\textbf{X})  &=& \frac{\sigma^2}{n^2h^{2d}(2\pi)^{2d}}\sum_{k=1}^n \Big| \int_{\mathbb{R}^d} e^{-i\langle \textbf{s},(\textbf{x}^*-\textbf{X}_k)\rangle /h}\frac{\Phi_K(\textbf{s}) }{\Phi_{\psi}(\frac{\textbf{s}}{h})} d\textbf{s} \Big|^2\frac{1}{\max\{f(\textbf{X}_k),f(\frac{1}{\bf{a}_n})\}^2}
\end{eqnarray*}
which yields
\begin{eqnarray*}
E[V(\hat{\theta}(\textbf{x}^*)|\textbf{X})] &=&\frac{\sigma^2}{nh^{d}(2\pi)^{2d}}  \int_{\mathbb{R}^d}\Big| \int_{\mathbb{R}^d} e^{-i\langle \textbf{s},(\textbf{x}^*/h-y)\rangle }\frac{\Phi_K(\textbf{s}) }{\Phi_{\psi}(\frac{\textbf{s}}{h})} d\textbf{s} \Big|^2\frac{f(h\textbf{y})}{\max\{f(h\textbf{\textbf{y}}),f(\frac{1}{\bf{a}_n})\}^2}d\textbf{y}.
\end{eqnarray*}
The variance of the conditional expectation is given by (observe again the definition of the weight $w_{n}$ in \eqref{wRD})
\begin{eqnarray*}
V(E[\hat{\theta}(\textbf{x}^*)|\textbf{X}]) &=& V\Big( \sum_{k=1}^n g(\textbf{X}_k) w_{n}(\textbf{x}^*,\textbf{X}_k)\Big)
\\ &=& \frac{1}{nh^{d}(2\pi)^{2d}}\int_{\mathbb{R}^d}\Big| \int_{\mathbb{R}^d} e^{-i\langle \textbf{s},(\textbf{x}^*/h-\textbf{y})\rangle }\frac{\Phi_K(\textbf{s}) }{\Phi_{\psi}(\frac{\textbf{s}}{h})} d\textbf{s} \Big|^2\frac{g^2(h\textbf{y})f(h\textbf{y})}{\max\{f(h\textbf{\textbf{y}}),f(\frac{1}{\bf{a}_n})\}^2}d\textbf{y}
\\ &-&\frac{1}{n(2\pi)^{2d}}\Big|\int_{\mathbb{R}^d} \int_{\mathbb{R}^d} e^{-i\langle \textbf{s},(\textbf{x}^*/h-\textbf{y})\rangle }\frac{\Phi_K(\textbf{s}) }{\Phi_{\psi}(\frac{\textbf{s}}{h})} d\textbf{s} \frac{g(h\textbf{y})f(h\textbf{y})}{\max\big\{f(h\textbf{\textbf{y}}),f(\frac{1}{\bf{a}_n})\big\}}d\textbf{y}\Big|^2,
\end{eqnarray*}
where the second summand is of order $O(n^{-1})$. Thus the variance can be written as
\begin{eqnarray}
 &&V(\hat{\theta}(\textbf{x}^*)) =E[V(\hat{\theta}(\textbf{x}^*)|\textbf{X})] + V(E[\hat{\theta}(\textbf{x}^*)|\textbf{X}]) \label{varberechnung}
\\ &=& \frac{1}{nh^{d}(2\pi)^{2d}}\int_{\mathbb{R}^d}\Big| \int_{\mathbb{R}^d} e^{-i\langle \textbf{s},(\textbf{x}^*/h-y)\rangle }\frac{\Phi_K(\textbf{s}) }{\Phi_{\psi}(\frac{\textbf{s}}{h})} d\textbf{s} \Big|^2\frac{(\sigma^2+g^2(h\textbf{y}))f(h\textbf{y})}{\max\{f(h\textbf{\textbf{y}}),f(\frac{1}{\bf{a}_n})\}^2}d\textbf{y}+O(n^{-1}) \nonumber
\end{eqnarray}
and the rate of convergence has a lower bound given by
$$
V(\hat{\theta}(\textbf{x}^*))^{-1/2}  =  \Omega\big( n^{1/2}h^{\beta+d/2}f(\textbf{a}_\textbf{n}^{-1})^{1/2} \big),
$$
where the symbol $b_n=\Omega(c_n)$ means that there exists a constant $C$ and $n_0 \in \mathbb{N}$ such that for all $n \ge n_0$ we have $|b_n|\ge C|c_n|$.
The variance has a lower bound
\begin{eqnarray*}
V(\hat{\theta}(\textbf{x}^*)) &\ge& \frac{1}{nh^{d}(2\pi)^{2d}}\int_{([\frac{-1}{ha_n},\frac{1}{ha_n}]^d)}\Big| \int_{\mathbb{R}^d} e^{-i\langle \textbf{s},(\textbf{x}^*/h-y)\rangle }\frac{\Phi_K(\textbf{s}) }{\Phi_{\psi}(\frac{\textbf{s}}{h})} d\textbf{s} \Big|^2\frac{(\sigma^2+g^2(h\textbf{y}))f(h\textbf{y})}{f(h\textbf{\textbf{y}})^2}d\textbf{y} \nonumber
\\ &\ge&  \frac{C}{nh^{d}(2\pi)^{2d}}\int_{([\frac{-1}{h a_n},\frac{1}{h a_n}]^d)}\Big| \int_{\mathbb{R}^d} e^{-i\langle \textbf{s},(\textbf{x}^*/h-y)\rangle }\frac{\Phi_K(\textbf{s}) }{\Phi_{\psi}(\frac{\textbf{s}}{h})} d\textbf{s} \Big|^2 d\textbf{y}  = C (nh^{d+2\beta})^{-1}(1+o(1)), \nonumber
\end{eqnarray*}
where we used Assumption \ref{assrv} and Parsevals equality. This yields to the upper bound
\begin{eqnarray}
V(\hat{\theta}(\textbf{x}^*))^{-1/2} &=& O\left(n^{1/2} h^{\beta+d/2}\right) \label{up}
\end{eqnarray}
For the proof of asymptotic normality we now show that the $l$-th cumulant of $G_l= \big| cum_l \big(V(\hat{\theta}(x^*))^{-1/2}\hat{\theta}(\textbf{x}^*)\big) \big|$  is vanishing asymptotically, whenever $l \ge 3$. For this purpose we recall the definition of the weights $w_{n}$ in \eqref{wRD} and obtain from \eqref{up} the estimate\allowdisplaybreaks{
\begin{eqnarray}
G_l &\le&C n^{l/2}h^{l\beta+dl/2} \sum_{k_1,...,k_l=1}^n |cum \big(Y_{k_1}w_{n}(\textbf{x}^*,\textbf{X}_{k_1}),...,Y_{k_l}w_{n}(\textbf{x}^*,\textbf{X}_{k_l})\big)| \label{52}
\\ &=& C n^{l/2}h^{l\beta+dl/2} \sum_{k=1}^n |cum_l \big(Y_{k}w_{n}(\textbf{x}^*,\textbf{X}_{k})\big)| \nonumber
\\ &=& C n^{l/2+1}h^{l\beta+dl/2} \sum_{\textbf{j}\in \{0,1\}^l} |( cum(U^{j_1}w_{n}(\textbf{x}^*,\textbf{X}_1),...,U^{j_l}w_{n}(\textbf{x}^*,\textbf{X}_1)) |, \nonumber
\end{eqnarray}
where we used (B2) and the notation $U^0=g(\textbf{X}_1)$ and $U^1=\epsilon$. This term can be written as
\begin{eqnarray*}
Cn^{l/2+1}h^{l\beta+dl/2}\sum_{s=0}^l \dbinom{l}{s} \sum_{\substack{ \textbf{j} \in\{0,1\}^l \\ j_1+...+j_l=s}}\big| cum(U^{j_1}w_{n}(\textbf{x}^*,\textbf{X}_1),...,U^{j_l}w_{n}(\textbf{x}^*,\textbf{X}_1)) \big|.
\end{eqnarray*}
By using the product theorem for cumulants [see e.g. \cite{bril2001}], we obtain
\begin{eqnarray}
 Cn^{l/2+1}h^{l\beta+dl/2}\sum_{s=0}^l \dbinom{l}{s} \sum_{\substack{ \textbf{j} \in\{0,1\}^l \\ j_1+...+j_l=s}} \Big| \sum_\nu \prod_{k=1}^p cum(A_{ij}, ij\in \nu_k) \Big|, \label{cum1}
\end{eqnarray}
where the third sum is calculated over all indecomposable partitions $\nu=(\nu_1,...,\nu_p)$ of the table
\begin{table}[htb!]
\centering
\begin{tabular}{l c c}
$A_{i1}$ &  & $ A_{i2}$ \\
\vdots & & \vdots \\
$A_{i1}$ &  & $ A_{i2}$ \\
& $A_{ij}$& \\
& \vdots & \\
& $A_{ij}$&
\end{tabular}
\end{table}
\\ (here the first $s$ rows have two and the last $l-s$ rows have one column) and
\begin{eqnarray*}
A_{i1} &=& \epsilon \qquad \qquad \qquad \qquad  \;1\le i \le s
\\ A_{i2} &=& w_{n}(\textbf{x}^*,\textbf{X}_1)) \qquad \qquad 1\le i \le s
\\ A_{ij} &=& g(\textbf{X}_1)w_{n}(\textbf{x}^*,\textbf{X}_1)) \qquad s+1 \le i \le l.
\end{eqnarray*}
As $\epsilon$ is independent of \textbf{X} only those indecomposable partitions yield a non zero cumulant, which seperate all $\epsilon$'s from the other terms. This means that for a partition $\nu$ there are $m(\nu)$ sets $\nu_1,...,\nu_{m(\nu)}$ which include only $\epsilon's$ while $\nu_{m(\nu)+1},...,\nu_p$ contain only $w_{n}(\textbf{x}^*,\textbf{X})$'s and $g(\textbf{X})w_{n}(\textbf{x}^*,\textbf{X})$'s. Thus \eqref{cum1} can be written as
\begin{eqnarray}
Cn^{l/2+1}h^{l\beta+dl/2}\sum_{s=0}^l \dbinom{l}{s} \sum_{\substack{ \textbf{j} \in\{0,1\}^l \\ j_1+...+j_l=s}} \Big| \sum_\nu \prod_{k=1}^{m(\nu)} cum_{s_k}(\epsilon) \prod_{k=m(\nu)+1}^p cum(A_{ij}, ij\in \nu_{k}) \Big|  \label{cum4}
\end{eqnarray}
with
\begin{eqnarray*}
A_{ij} &=& w_{n}(\textbf{x}^*-\textbf{X}_1)) \qquad \qquad 1\le i \le s
\\ A_{ij} &=& g(\textbf{X})w_{n}(\textbf{x}^*-\textbf{X}_1)) \qquad s+1 \le i \le l.
\end{eqnarray*}
and $s_1+...+s_{m(\nu)}=s$ . Furthermore we have $s_i \ge 2$, because the noise terms $\epsilon$ have mean zero, and each set $\nu_{m(\nu)+1},...,\nu_p$ includes at least one $A_{ij}$ with $1\le i \le s$ because otherwise the partition would not be indecomposable. Let $a_{r}=|\nu_{r}|$ denote the number of elements in the set $\nu_r \mbox{ } (r=m(\nu)+1,...,p)$, then we get $a_{m+1}+...+a_p=l$. Furthermore for $r\in\{m+1,...,p\} \mbox{ the cumulant }  cum(A_{ij},ij\in \nu_r)$ equals
\begin{eqnarray}
 cum(g(\textbf{X}_1)w_{n}(\textbf{x}^*,\textbf{X}_1)),...,g(\textbf{X}_1)w_{n}(\textbf{x}^*,\textbf{X}_1)),w_{n}(\textbf{x}^*,\textbf{X}_1)),...,w_{n}(\textbf{x}^*,\textbf{X}_1))) \label{cum2}
\end{eqnarray}
because of the symmetry of the arguments in the cumulant. In the next step we denote by $b_r$ the number of components of the form $g(\textbf{X}_1)w_n(\textbf{x}^*,\textbf{X}_1)$ and show the estimate
\begin{eqnarray}
E\Big[\prod_{i=1}^{b_r}|g(\textbf{X}_1)w_{n}(\textbf{x}^*,\textbf{X}_1))|\prod_{j=1}^{a_r-b_r}|w_{n}(\textbf{x}^*,\textbf{X}_1))|\Big] \le  \frac{C^{a_r}}{n^{a_r}h^{a_r(\beta+d)}f(\frac{1}{\bf{a}_n})^{a_r}} \label{cum3}
\end{eqnarray}
(which does not depend on $b_r$). From Theorem 5.1 we then obtain that the term in \eqref{cum2} is of order $O(n^{-a_r}h^{-a_r(\beta+d)}f(1/\textbf{a}_\textbf{n})^{-a_r})$. Equations \eqref{52}, \eqref{cum4} and \eqref{cum2} yield for the cumulants of order $l\ge 3$
\begin{eqnarray*}
G_l &\le& Cn^{l/2+1}h^{l\beta+dl/2}\sum_{s=0}^l \dbinom{l}{s} \sum_{\substack{ \textbf{j} \in\{0,1\}^l \\ j_1+...+j_l=s}} \Big| \sum_\nu \prod_{k=1}^{m(\nu)} cum_{s_k}(\epsilon)\prod_{r=m(\nu)+1}^p \frac{C^{a_r}}{n^{a_r}h^{a_r(d+\beta)}f(\frac{1}{\bf{a}_n})^{a_r}}\Big|
\\ &=& O\big((n^{l/2-1}h^{ld/2}f({\bf{a}_n}^{-1})^l)^{-1}\big)=o(1),
\end{eqnarray*}
which shows the asymptotic normality.
\\ In order to prove the remaining estimate \eqref{cum3} we use the definition of $w_{n}(\textbf{x}^*,\textbf{X}_1)$ and obtain for the term on the left hand side of \eqref{cum3}
\begin{eqnarray*}
L_n &=& \int_{\mathbb{R}^d} |g(\textbf{x})|^{b_r} \Big(\frac{1}{nh^d(2\pi)^d} \Big|\int_{\mathbb{R}^d} e^{-i\langle \textbf{s},(\textbf{x}^*-\textbf{x})\rangle /h}\frac{\Phi_K(\textbf{s}) }{\Phi_{\psi}(\frac{\textbf{s}}{h})}\frac{1}{\max\{f(\textbf{x}),f(\frac{1}{\bf{a}_n})\}} d\textbf{s}\Big|\Big)^{a_r}f(\textbf{x}) d\textbf{x}
\\ &\le&\frac{C}{n^{a_r}h^{a_rd}} \int_{\mathbb{R}^d} |g(\textbf{x})|^{b_r} \Big( \int_{\mathbb{R}^d}\frac{|\Phi_K(\textbf{s})| }{|\Phi_{\psi}(\frac{\textbf{s}}{h})|}\frac{1}{\max\{f(\textbf{x}),f(\frac{1}{\bf{a}_n})\}} d\textbf{s}\Big)^{a_r}f(\textbf{x}) d\textbf{x}
\\ &\le& \frac{C}{n^{a_r}h^{a_rd}f(\frac{1}{\bf{a}_n})^{a_r}}\int_{\mathbb{R}^d}  \Big( \int_{\mathbb{R}^d}\frac{|\Phi_K(\textbf{s})| }{|\Phi_{\psi}(\frac{\textbf{s}}{h})|} d\textbf{s}\Big)^{a_r}f(\textbf{x}) d\textbf{x},
\end{eqnarray*}
where we used the fact that $g$ is bounded. Using this inequality and Assumption \ref{asspsiadd}(A) it follows that $L_n \le C/n^{a_r}h^{a_r(d+\beta)}f(\frac{1}{\bf{a}_n})^{a_r},$ which proves \eqref{cum3}.
}

\medskip

\textbf{Proof of Lemma \ref{lemaddrd}:}
Similar to the proof of Theorem \ref{theoaddrd}, we have to calculate the cumulants of the estimators $\hat{\alpha}_{j,Q_{I_j^c}}(\textbf{x}_{I_j}^*)$. We start with the first order cumulant
\begin{eqnarray*}
E[\hat{\alpha}_{j,Q_{I_j^c}}(\textbf{x}_{I_j}^*)] &=& \frac{1}{h^d(2\pi)^d}\int_{\mathbb{R}^{d-d_j}} \int_{\mathbb{R}^d} \int_{\mathbb{R}^d} g(\textbf{x}) \int_{\mathbb{R}^d} \frac{\Phi_K(\textbf{s}) }{\Phi_{\psi}(\frac{\textbf{s}}{h})}\frac{e^{-i\langle \textbf{s},(\textbf{x}^*-\textbf{x})\rangle /h}f(\textbf{x})}{\max\{f(\textbf{x}),f(\frac{1}{\bf{a}_n})\}} d\textbf{s}d\textbf{x}dQ_{I_j^c}(\textbf{x}_{I_j^c}^*)
\end{eqnarray*}
and  with the same arguments as in the proof of Theorem \ref{theoaddrd},  we obtain a bias of order $O(h^{s-1})$. For the calculation of the variance of $\hat{\alpha}_{j,Q_{I_j^c}}(\textbf{x}_{I_j}^*)$ we investigate its conditional variance. Recalling the definitions \eqref{wRD} and \eqref{LI} it follows by a straightforward argument
\begin{eqnarray*}
V(\hat{\alpha}_{j,Q_{I_j^c}}(\textbf{x}_{I_j}^*)|\textbf{X}) &=&  \frac{\sigma^2}{n^2h^{2d}(2\pi)^{2d}}\sum_{k=1}^n  \Big|\int_{\mathbb{R}^d} e^{-i\langle \textbf{w},\textbf{X}_k\rangle /h} e^{i\langle \textbf{w}_{I_j},\textbf{x}_{I_j}\rangle /h}L_{I_j^c}\Big(\frac{\textbf{w}_{I_j^c}}{h}\Big)\frac{\Phi_K(\textbf{w})}{\Phi_\psi(\frac{\textbf{w}}{h})} d\textbf{w}\Big|^2
\\ && \times \frac{1}{\max\{f(\textbf{X}_{k}),f(\frac{1}{\bf{a}_n})\}^2},
\end{eqnarray*}
which gives
\begin{eqnarray*}
E\big[V(\hat{\alpha}_{j,Q_{I_j^c}}(\textbf{x}_{I_j}^*)|\textbf{X})\big] &=&  \frac{\sigma^2}{nh^{2d}(2\pi)^{2d}}\int_{\mathbb{R}^d} \Big|\int_{\mathbb{R}^d} e^{-i\langle \textbf{w},\textbf{x}\rangle /h} e^{i\langle \textbf{w}_{I_j},\textbf{x}_{I_j}\rangle /h}L_{I_j^c}\Big(\frac{\textbf{w}_{I_j^c}}{h}\Big)\frac{\Phi_K(\textbf{w})}{\Phi_\psi(\frac{\textbf{w}}{h})} d\textbf{w}\Big|^2
\\ && \times \frac{f(\textbf{x})}{\max\{f(\textbf{x}),f(\frac{1}{\bf{a}_n})\}^2} d\textbf{x}.
\end{eqnarray*}
The variance of the conditional expectation can be calculated as
\begin{eqnarray*}
V\big(E[\hat{\alpha}_{j,Q_{I_j^c}}(\textbf{x}_{I_j}^*)|\textbf{X}]\big) &=& \frac{1}{nh^{d}(2\pi)^{2d}}\int_{\mathbb{R}^d} \Big|\int_{\mathbb{R}^d} e^{-i\langle \textbf{w},\textbf{x}\rangle } e^{i\langle \textbf{w}_{I_j},\textbf{x}_{I_j}\rangle /h}L_{I_j^c}\Big(\frac{\textbf{w}_{I_j^c}}{h}\Big)\frac{\Phi_K(\textbf{w})}{\Phi_\psi(\frac{\textbf{w}}{h})} d\textbf{w}\Big|^2
\\ && \times \frac{g(h\textbf{x})^2f(h\textbf{x})}{\max\{f(h\textbf{x}),f(\frac{1}{\bf{a}_n})\}^2} d\textbf{x}
\\ &-& \frac{1}{n(2\pi)^{2d}}\Big|\int_{\mathbb{R}^d} \int_{\mathbb{R}^d} e^{-i\langle \textbf{w},\textbf{x}\rangle } e^{i\langle \textbf{w}_{I_j},\textbf{x}_{I_j}\rangle /h}L_{I_j^c}\Big(\frac{\textbf{w}_{I_j^c}}{h}\Big)\frac{\Phi_K(\textbf{w})}{\Phi_\psi(\frac{\textbf{w}}{h})} d\textbf{w}
\\ && \times \frac{g(h\textbf{x})f(h\textbf{x})}{\max\{f(h\textbf{x}),f(\frac{1}{\bf{a}_n})\}} d\textbf{x}\Big|^2,
\end{eqnarray*}
where the second summand is of order $O(n^{-1})$. Therefore it follows
\begin{eqnarray*}
V(\hat{\alpha}_{j,Q_{I_j^c}}(\textbf{x}_{I_j}^*))  &=& \frac{1}{nh^{d}(2\pi)^{2d}}\int_{\mathbb{R}^d} \Big|\int_{\mathbb{R}^d} e^{-i\langle \textbf{w},\textbf{y}\rangle } e^{i\langle \textbf{w}_{I_j},\textbf{x}_{I_j}\rangle /h}L_{I_j^c}\Big(\frac{\textbf{w}_{I_j^c}}{h}\Big)\frac{\Phi_K(\textbf{w})}{\Phi_\psi(\frac{\textbf{w}}{h})} d\textbf{w}\Big|^2
\\ && \times \frac{(\sigma^2+g(h\textbf{y})^2)f(h\textbf{y})}{\max\{f(h\textbf{y}),f(\frac{1}{\bf{a}_n})\}^2} d\textbf{y}+O(n^{-1}).
\end{eqnarray*}
The upper bound for this term is obtained from Assumption \ref{assrv} which gives
\begin{eqnarray}
\frac{(\sigma^2+g(h\textbf{y})^2)f(h\textbf{y})}{\max\{f(h\textbf{y}),f(\frac{1}{\bf{a}_n})\}^2}=O({\bf a_n}^{-1}). \label{absch}
\end{eqnarray}
Therefore an application of Parseval's equality and Assumption \ref{asspsi}(C) yields
\begin{eqnarray}
 V(\hat{\alpha}_{j,Q_{I_j^c}}(\textbf{x}_{I_j}^*) \le \frac{C}{nh^{d+2\beta-\gamma_j}f(\frac{1}{\bf{a}_n})}. \label{wn}
\end{eqnarray}
A similar argument as in the proof of Theorem \ref{theoaddrd} gives the lower bound $V(\hat{\alpha}_{j,Q_{I_j^c}}(\textbf{x}_{I_j}^*) \ge C/nh^{d+2\beta-\gamma_j}.$ Finally the statement that the $l$-th cumulant of $V(\hat{\alpha}_{j,Q_{I_j^c}}(\textbf{x}_{I_j}^*)^{-1/2}\hat{\alpha}_{j,Q_{I_j^c}}(x_{I_j}^*)$ is of order $o(1)$ can be shown by similiar arguments as in the proof of Theorem \ref{theoaddrd}.

\medskip

\textbf{Proof of Theorem \ref{theoremaddrd}: }
The proof follows by similar arguments as given in the previous Sections. For the sake of brevity we restrict ourselves for the calculation of the first and second order cumulants. For this purpose we show, that the estimate $\hat{c}$ has a faster rate of convergence than $\hat{\alpha}_{j,Q_{I_j^c}}(\textbf{x}_{I_j}^*)$ for at least one $j \in \{1,...,m\}$. If this statement is correct the asymptotic variance of the statistic
\begin{eqnarray*}
\hat{\theta}^{add,RD}(\textbf{x}^*) &=& \sum_{j=1}^m \hat{\alpha}_{j,Q_{I_j^c}}(\textbf{x}_{I_j}^*) -(m-1)\hat{c}
\end{eqnarray*}
is determined by its first term. Recalling the notation \eqref{thetaddrd} this term has the representation
\begin{eqnarray}
\hat{D}_n = \sum_{j=1}^m \hat{\alpha}_{j,Q_{I_j^c}}(\textbf{x}_{I_j}^*) &=&\sum_{j=1}^m \sum_{k=1}^n Y_k  w_{n}^{add,RD}(\textbf{x}_{I_j}^*,\textbf{X}_k) \label{Dn}
\end{eqnarray}
and can be treated in the same way as before. The resulting bias of $\hat{D}_n$ is the sum of the biases of the individual term and therefore also of order $O(h^{s-1})$. The conditional variance is given by
\begin{eqnarray*}
V(\hat{D}_n|\textbf{X}) &=& \sigma^2 \sum_{k=1}^n \Big|\sum_{j=1}^m w_{n}^{add,RD}(\textbf{x}_{I_j}^*,\textbf{X}_k)\Big|^2
\\ &=& \frac{\sigma^2}{n^2h^{2d}(2\pi)^{2d}}\sum_{k=1}^n \Big|  \int_{\mathbb{R}^d} e^{i\langle \textbf{w},\textbf{X}_k\rangle /h} \Big( \sum_{j=1}^m e^{-i\langle \textbf{w}_{I_j},\textbf{x}_{I_j}^*\rangle /h}L_{I_j^c}\Big(\frac{\textbf{w}_{I_j^c}}{h} \Big)\Big) \frac{\Phi_K(\textbf{w})}{\Phi_\psi(\frac{\textbf{w}}{h})} d\textbf{w}
\\ && \times\frac{1}{\max\{ f(\textbf{X}_k),f(\frac{1}{\bf{a}_n})\}} \Big|^2.
\end{eqnarray*}
This yields for expectation of the conditional variance
\begin{eqnarray*}
\lefteqn{E\big[V(\hat{D}_n|\textbf{X})\big]}
\\ &=& \frac{\sigma^2}{nh^{2d}(2\pi)^{2d}}\int_{\mathbb{R}^d} \Big|  \int_{\mathbb{R}^d} e^{i\langle \textbf{w},\textbf{s}\rangle /h} \Big( \sum_{j=1}^m e^{-i\langle \textbf{w}_{I_j},\textbf{x}_{I_j}^*\rangle /h}L_{I_j^c}\Big(\frac{\textbf{w}_{I_j^c}}{h} \Big)\Big) \frac{\Phi_K(\textbf{w})}{\Phi_\psi(\frac{\textbf{w}}{h})} d\textbf{w} \Big|^2 \frac{f(\textbf{s})}{\max\{ f(\textbf{s}),f(\frac{1}{\bf{a}_n})\}^2}d\textbf{s}
\end{eqnarray*}
and the variance of the conditional expectation is obtained as
\begin{eqnarray*}
\lefteqn{V\big(E[\hat{D}_n|\textbf{X}]\big)}
\\ &=& \frac{1}{nh^{d}(2\pi)^{2d}}\int_{\mathbb{R}^d} \Big|  \int_{\mathbb{R}^d} e^{i\langle \textbf{w},\textbf{s}\rangle } \Big( \sum_{j=1}^m e^{-i\langle \textbf{w}_{I_j},\textbf{x}_{I_j}^*\rangle /h}L_{I_j^c}\Big(\frac{\textbf{w}_{I_j^c}}{h} \Big)\Big) \frac{\Phi_K(\textbf{w})}{\Phi_\psi(\frac{\textbf{w}}{h})} d\textbf{w} \Big|^2 \frac{g(h\textbf{s})^2f(h\textbf{s})}{\max\{ f(h\textbf{s}),f(\frac{1}{\bf{a}_n})\}^2}d\textbf{s}
\\ &-& \frac{1}{n(2\pi)^{2d}} \Big|\int_{\mathbb{R}^d}   \int_{\mathbb{R}^d} e^{i\langle \textbf{w},\textbf{s}\rangle } \Big( \sum_{j=1}^m e^{-i\langle \textbf{w}_{I_j},\textbf{x}_{I_j}^*\rangle /h}L_{I_j^c}\Big(\frac{\textbf{w}_{I_j^c}}{h} \Big)\Big) \frac{\Phi_K(\textbf{w})}{\Phi_\psi(\frac{\textbf{w}}{h})} d\textbf{w}  \frac{g(h\textbf{s})f(h\textbf{s})}{\max\{ f(h\textbf{s}),f(\frac{1}{\bf{a}_n})\}}d\textbf{s}\Big|^2 ,
\end{eqnarray*}
where the second summand is of order $O(n^{-1})$. Thus yields for the variance
\begin{eqnarray*}
V(\hat{D}_n) &=& \frac{1}{nh^{2d}(2\pi)^{2d}}\int_{\mathbb{R}^d} \Big|  \int_{\mathbb{R}^d} e^{i\langle \textbf{w},\textbf{s}\rangle /h} \Big( \sum_{j=1}^m e^{-i\langle \textbf{w}_{I_j},\textbf{x}_{I_j}^*\rangle /h}L_{I_j^c}\Big(\frac{\textbf{w}_{I_j^c}}{h} \Big)\Big) \frac{\Phi_K(\textbf{w})}{\Phi_\psi(\frac{\textbf{w}}{h})} d\textbf{w} \Big|^2
\\ && \times \frac{(\sigma^2+g(\textbf{s})^2)f(\textbf{s})}{\max\{ f(\textbf{s}),f(\frac{1}{\bf{a}_n})\}}d\textbf{s}+O(n^{-1})
\end{eqnarray*}
In order to obtain bounds for the rate of the variance, we use the lower bound for $\max\{ f(h\textbf{s}),f(\frac{1}{a_n})\}$ mentioned in \eqref{absch} and Parseval's equality which yields
\begin{eqnarray*}
 \Big(\frac{1}{nh^{d}f(\frac{1}{\bf{a}_n})}  \int_{\mathbb{R}^d} \left| \sum_{j=1}^m e^{-i\langle \textbf{w}_{I_j},\textbf{x}_{I_j}^*\rangle /h}L_{I_j^c}\Big(\frac{\textbf{w}_{I_j^c}}{h} \Big)\right|^2 \frac{|\Phi_K(\textbf{w})|^2}{|\Phi_\psi(\frac{\textbf{w}}{h})|^2} d\textbf{w}\Big)^{1/2} &=& O\left( (nh^{d+2\beta-\gamma_{min}}f(\bf{a}_n^{-1}))^{-1}\right)
\end{eqnarray*}
as an upper bound, where the last estimate follows from Assumption \ref{asspsi}. The lower bound is of order $\Omega( (nh^{d+2\beta-\gamma_{min}})^{-1})$, where we use Assumption \ref{assrv} and the same calculations as in the previous Section.
These are in fact the same bounds as for $\hat{\alpha}_{j^*,Q_{I_{j^*}^c}}(\textbf{x}_{I_{j^*}}^*)$ with $j^* = argmin_j \:\gamma_j$. This means that
\[ \hat{D}_n-E[\hat{D}_n] = O_P(n^{-1/2}h^{-d/2-\beta-\gamma_{min}/2}f({\bf a_n}^{-1})^{-1/2}) \]
 In the last step we show that the estimate $\hat{c}$ has a faster rate of convergence. For this purpose we write  $\hat{c}$ as weighted sum of independent random variables that is
\begin{eqnarray*}
\hat{c} &=& \int_{\mathbb{R}^d} \hat{\theta}(\textbf{x}^*) dQ(\textbf{x}^*)= \frac{1}{nh^d(2\pi)^d} \sum_{k=1}^n   \int_{\mathbb{R}^d} e^{i\langle \textbf{w},\textbf{X}_k\rangle /h}\Big(\prod_{j=1}^m L_{I_j^c}\Big( \frac{\textbf{w}_{I_j^c}}{h}\Big)\Big)\frac{\Phi_K(\textbf{w})}{\Phi_\psi(\frac{\textbf{w}}{h})} d\textbf{w} \frac{Y_k}{\max\{ f(\textbf{X}_k),f(\frac{1}{\bf{a}_n})\}}
\end{eqnarray*}
It now follows by similar calculations as given in the previous paragraph and Assumption \ref{asspsi}(C) that
\[ V(\hat{c})=o(V(\sum_{j=1}^m \hat{\alpha}_{j,Q_{I_j^c}}(\textbf{x}_{I_j}^*)))\]
and thus we can ignore  the term $\hat{c}$ for the calculation of the asymptotic variance of the statistic $\hat{\theta}^{add,RD}$.

\medskip
\textbf{Proof of Lemma \ref{lemfd}: }
\allowdisplaybreaks{
Observing the representation \eqref{zmarg} and \eqref{sumFD} we decompose the estimator into its deterministic and stochastic part, that is
\begin{eqnarray}
 \label{dspart}
\hat{\theta}_{I_j}^{FD}(\textbf{x}_{I_j}^*) &=& \hat{E}_{1n}+\hat{E}_{2n}
\end{eqnarray}
where
\begin{eqnarray*}
\hat{E}_{1n} &=&  \frac{1}{(2n+1)^{d-d_j}} \sum_{\textbf{k} \in \{-n,...,n\}^d} (g_{I_1}(\textbf{z}_{\textbf{k}_{I_1}})+...+g_{I_m}(\textbf{z}_{\textbf{k}_{I_m}})) w_{\textbf{k}_{I_j},n}(\textbf{x}_{I_j}^*) \nonumber
\\  \hat{E}_{2n} &=& \frac{1}{(2n+1)^{d-d_j}}\sum_{\textbf{k} \in \{-n,...,n\}^d}^n \epsilon_{\textbf{k}} w_{\textbf{k}_{I_j},n}(\textbf{x}_{I_j}^*) \nonumber
\end{eqnarray*}
and $w_{\textbf{k}_{I_j},n}(\textbf{x}_{I_j}^*)$ are defined in \eqref{wfdadd}. In a first step we show, that the bias of $\hat{\theta}_{I_j}^{FD}$ is of order $O( \frac{1}{n^2h^{2+d_j+\beta_j}a_n^3})$. For this purpose we rewrite the deterministic part as
\begin{eqnarray*}
\hat{E}_{1n} &=& \hat{E}_{1n}^{(1)}+\hat{E}_{1n}^{(2)}
\end{eqnarray*}
where
\begin{eqnarray*}
\hat{E}_{1n}^{(1)} &=&\sum_{\textbf{k}_{I_j} \in \{-n,...,n\}^{d_j}} g_{I_j}(\textbf{z}_{\textbf{k}_{I_j}})w_{\textbf{k}_{I_j},n}(\textbf{x}_{I_j})
\\\hat{E}_{1n}^{(2)} &=& \frac{1}{(2n+1)^{d-d_j}}\sum_{\textbf{k}_{I_j^c} \in \{-n,...,n\}^{d-d_j}}  \left(g_{I_1}(\textbf{z}_{\textbf{k}_{I_1}})+...+g_{I_{j-1}}(\textbf{z}_{\textbf{k}_{I_{j-1}}})+g_{I_{j+1}}(\textbf{z}_{\textbf{k}_{I_{j+1}}})+...+g_{I_m}(\textbf{z}_{\textbf{k}_{I_m}})\right)
\\ &&\times \sum_{\textbf{k}_{I_j} \in \{-n,...,n\}^{d_j}} w_{\textbf{k}_{I_j},n}(\textbf{x}_{I_j}^*),
\end{eqnarray*}
where the second summand is of order $\hat{E}_{1n}^{(2)} = o(\frac{a_n^{r-d_j}}{h^{\beta_j+d_j}})=O(h^{s-1})$, which follows from Assumption \ref{asstheta}(D).
For the difference of the first summand and $\theta_{I_1}^{add}(\textbf{x}_{I_j}^*)$ we use the same calculation as in \cite{bb2008} and obtain
\begin{eqnarray*}
\hat{E}_{1n}^{(1)}-\theta_{I_j}^{add}(\textbf{x}_{I_j}) &=& O(h^{s-1}) +O(\frac{1}{n^2a_n^2h^{d_j+2+\beta_j}}).
\end{eqnarray*}
Note that the Rieman-approximation does not provide an error of order $O((na_n)^{-d})$, but we can show that the lattice structure yields an error term of order $O( (n^2h^2a_n^3)^{-1})$.
In the next step we derive the variance of the estimator $\hat{\theta}_{I_j}^{FD}$. We can neglect the deterministic part $\hat{E}_{2n}$ in \eqref{dspart} and obtain from Parseval's equality and Assumption \ref{asspsiadd}(B)
\begin{eqnarray}
V(\hat{\theta}_{I_j}^{FD}(\textbf{x}_{I_j}^*)) &=&  \frac{\sigma^2}{(2n+1)^{d-d_j}} \sum_{\textbf{k}_{I_j} \in \{-n,...,n\}^{d_j}} |w_{\textbf{k}_{I_j},n}(\textbf{x}_{I_j}^*)|^2 \nonumber
\\ &=& \frac{\sigma^2}{(2n+1)^{d-d_j}n^{2d_j}h^{2d_j}a_n^{2d_j}(2\pi)^{2d_j}}\sum_{\textbf{k}_{I_j} \in \{-n,...,n\}^{d_j}} \Big| \int_{\mathbb{R}^{d_j}} e^{-i\langle \textbf{w}, (\textbf{x}_{I_j}^*-\textbf{z}_{k_{I_j}})\rangle /h}  \frac{\Phi_K(\textbf{w}) }{\Phi_{\psi_{I_j}}(\frac{\textbf{w}}{h})} d\textbf{w} \Big|^2 \nonumber
\\ &=&  \frac{\sigma^2}{(2n+1)^{d-d_j}n^{d_j}h^{d_j}a_n^{d_j}(2\pi)^{2d_j}} \nonumber
\\ && \Big(\int_{[-1/(ha_n),1/(ha_n)]^{d_j}} \Big| \int_{\mathbb{R}^{d_j}} e^{-i\langle \textbf{w}, (\textbf{x}_{I_j}^*/h-\textbf{s})\rangle }  \frac{\Phi_K(\textbf{w}) }{\Phi_{\psi_{I_j}}(\frac{\textbf{w}}{h})} d\textbf{w} \Big|^2 ds \nonumber +O((na_n)^{-1})\Big) \nonumber
\\ &\sim& \frac{\sigma^2}{(2n+1)^{d-d_j}n^{d_j}h^{d_j}a_n^{d_j}(2\pi)^{2d_j}}\int_{\mathbb{R}^{d_j}} \Big| \int_{\mathbb{R}^{d_j}} e^{-i\langle \textbf{w}, (\textbf{x}_{I_j}^*/h-\textbf{s})\rangle }  \frac{\Phi_K(\textbf{w}) }{\Phi_{\psi_{I_j}}(\frac{\textbf{w}}{h})} d\textbf{w} \Big|^2 ds (1+o(1)) \nonumber
\\ &=& \frac{\sigma^2}{(2n+1)^{d-d_j}n^{d_j}h^{d_j}a_n^{d_j}(2\pi)^{2d_j}}\int_{\mathbb{R}^{d_j}}   \frac{|\Phi_K(\textbf{w})|^2 }{|\Phi_{\psi_{I_j}}(\frac{\textbf{w}}{h})|^2} d\textbf{w}  (1+o(1)) \nonumber
\\ &=& \frac{\sigma^2C}{(2n+1)^{d-d_j}n^{d_j}h^{d_j+2\beta_j}a_n^{d_j}(2\pi)^{2d_j}} \sim \frac{C}{n^dh^{d_j+2\beta_j}a_n^{d_j}}. \label{vtheta} \nonumber
\end{eqnarray}
 For the proof of the asymptotic normality, we finally show that the $l$-th cumulant of $V(\hat{\theta}_{I_j}^{FD}(\textbf{x}_{I_j}^*))^{-1/2} \hat{\theta}_{I_j}^{FD}(\textbf{x}_{I_j}^*)$ converges to zero for $l \ge 3$, which  completes the proof of Lemma \ref{lemfd}. For this purpose we note that

\begin{eqnarray*}
&& |cum_l(V(\hat{\theta}_{I_j}^{FD}(\textbf{x}_{I_j}^*))^{-1/2} \hat{\theta}_{I_j}^{FD}(\textbf{x}_{I_j}^*))| \le |C n^{ld/2}h^{ld_j/2+l\beta_j}a_n^{ld_j/2} cum_l(\hat{\theta}_{I_j}^{FD}(\textbf{x}_{I_j}^*))|
\\ &\le& \Big|\frac{C}{(n^{ld/2}h^{ld_j/2}a_n^{ld_j/2}}  \sum_{\textbf{k}_1,...,\textbf{k}_l \in \{-n,...,n\}^d} \prod_{m=1}^l \left( \int_{\mathbb{R}^{d_j}}  e^{-i\langle \textbf{w},(\textbf{x}_{I_j}^*-z_{k_{m,I_j}})\rangle /h} \frac{\Phi_K(\textbf{w})}{\Phi_\psi(\frac{\textbf{w}}{h})} d\textbf{w}\right) cum(\epsilon_{\textbf{k}_1},...,\epsilon_{\textbf{k}_l}) \Big|
\\ &\le& \frac{C}{n^{ld/2}h^{ld_j/2}a_n^{ld_j/2}}  \sum_{\textbf{k}_1 \in \{-n,...,n\}^d}
 \prod_{m=1}^l \left( \int_{\mathbb{R}^{d_j}}   \frac{|\Phi_K(\textbf{w})|}{|\Phi_\psi(\frac{\textbf{w}}{h})|} d\textbf{w}\right),
\end{eqnarray*}
where $\kappa_l$ denotes the $l$-th cumulant of $\epsilon$. From Assumption 1 it follows that this term is bounded by
\[  \frac{C }{n^{ld/2}h^{ld_j/2}a_n^{ld_j}/2}(2n+1)^d h^{-l\beta_j} = C n^{-ld/2+1}h^{-ld_j/2}a_n^{ld_j/2}, \]
which converges to zero for $l \ge 3$.

\medskip

\textbf{Proof of Theorem \ref{theofd}: }
In the following discussion we ignore the constant term $g_0=\theta_0$ because the mean
\[ \hat{\theta}_0 =\frac{1}{n^d} \sum_{\textbf{k} \in \{-n,...,n\}^d} Y_{\textbf{k}} \]
is a $\sqrt{n^d}$-consistent estimator for this constant and the nonparametric components in \eqref{fixadd} can only be estimated
at slower rates. Note that
\begin{eqnarray*}
\hat{\theta}^{add,FD}(\textbf{x}^*)  &=&  \sum_{\textbf{k} \in \{-n,...,n\}^d} Y_\textbf{k} \sum_{j=1}^n \frac{1}{(2n+1)^{d-d_j}}w_{\textbf{k}_{I_j},n}(\textbf{x}_{I_j}^*)
\end{eqnarray*}
and obtain the asymptotic distribution with the same arguments as in the proof of Lemma \ref{lemfd}.
}

\medskip

\textbf{Proof of Theorem \ref{lemfd2}: }
Under the assumption of an MA(q)-dependency structure \eqref{MAq} there are no changes in the calculation of the mean of the estimator $\hat{\theta}_{I_j}^{FD}$ and we only have to calculate the cumulants of order $l\ge2$ in order to establish the asymptotic normality. We start with the variance, which is given by
\begin{eqnarray*}
V(\hat{\theta}_{I_j}^{FD}(\textbf{x}_{I_j}^*)) &=& \frac{1}{(2n+1)^{2(d-d_j)} }\sum_{\textbf{k}_1,\textbf{k}_2 \in \{-n,...,n\}^d} w_{\textbf{k}_{1,I_j},n}(\textbf{x}_{I_j}^*)\overline{w}_{\textbf{k}_{2,I_j},n}(\textbf{x}_{I_j}^*) cum(\epsilon_{\textbf{k}_1},\epsilon_{\textbf{k}_2})
\\ &=& \frac{1}{(2n+1)^{2(d-d_j)} }\sum_{\textbf{k}_1 \in \{-n,...,n\}^d} \sum_{\textbf{k}_2 : \parallel \textbf{k}_2-\textbf{k}_1 \parallel_\infty \le 2q} \sum_{\textbf{r}_1,\textbf{r}_2 \in \{-q,...,q\}^d} w_{\textbf{k}_{1,I_j},n}(\textbf{x}_{I_j}^*)\overline{w}_{\textbf{k}_{2,I_j},n}(\textbf{x}_{I_j}^*)
\\ && cum(\beta_{\textbf{r}_1} Z_{\textbf{k}_1-\textbf{r}_1}, \beta_{\textbf{r}_2} Z_{\textbf{k}_2-\textbf{r}_2})
\\ &=& \frac{1}{(2n+1)^{2(d-d_j)} }\sum_{\textbf{k}_1 \in \{-n,...,n\}^d} \sum_{\textbf{k}_2 : \parallel \textbf{k}_2-\textbf{k}_1 \parallel_\infty \le 2q} \sum_{\textbf{r}_1 \in \{-q,...,q\}^d} w_{\textbf{k}_{1,I_j},n}(\textbf{x}_{I_j}^*)\overline{w}_{\textbf{k}_{2,I_j},n}(\textbf{x}_{I_j}^*)
\\ && cum(\beta_{\textbf{r}_1} Z_{\textbf{k}_1-\textbf{r}_1}, \beta_{\textbf{k}_2-\textbf{k}_1+\textbf{r}_1} Z_{\textbf{k}_1-\textbf{r}_1})
\\ &=& \frac{\sigma^2}{(2n+1)^{2(d-d_j)} }\sum_{\textbf{k}_1 \in \{-n,...,n\}^d} \sum_{\textbf{k}_2 : \parallel \textbf{k}_2-\textbf{k}_1 \parallel_\infty \le 2q} \sum_{\textbf{r}_1 \in \{-q,...,q\}^d} w_{\textbf{k}_{1,I_j},n}(\textbf{x}_{I_j}^*)\overline{w}_{\textbf{k}_{2,I_j},n}(\textbf{x}_{I_j}^*)
\\ && \beta_{\textbf{r}_1} \beta_{\textbf{k}_2-\textbf{k}_1+\textbf{r}_1}
\\ &=&  \frac{\sigma^2}{(2n+1)^{2(d-d_j)} }\sum_{\textbf{k}_1 \in \{-n,...,n\}^d} \sum_{\substack{ \textbf{l} \in \mathbb{Z}^d \\  \parallel \textbf{l} \parallel_\infty \le 2q}} \sum_{\textbf{r}_1 \in \{-q,...,q\}^d} w_{\textbf{k}_{1,I_j},n}(\textbf{x}_{I_j}^*)\overline{w}_{\textbf{l}_{I_j}+\textbf{k}_{1,I_j},n}(\textbf{x}_{I_j}^*)
\\ && \beta_{\textbf{r}_1} \beta_{\textbf{l}+\textbf{r}_1}
\\ &=& \frac{\sigma^2(1+o(1))}{(2n+1)^{(d-d_j)} }\sum_{\textbf{k}_{1,I_j} \in \{-n,...,n\}^{d_j}} \sum_{\substack{ \textbf{l} \in \mathbb{Z}^d \\  \parallel \textbf{l} \parallel_\infty \le 2q}} \sum_{\textbf{r}_1 \in \{-q,...,q\}^d} |w_{\textbf{k}_{1,I_j},n}(\textbf{x}_{I_j}^*)|^2 \beta_{\textbf{r}_1} \beta_{\textbf{l}+\textbf{r}_1},
\end{eqnarray*}
where we used a Taylor-approximation for the weights $ \overline{w}_{\textbf{l}_{I_j}+\textbf{k}_{1,I_j},n}(\textbf{x}_{I_j}^*)  = \overline{w}_{\textbf{k}_{1,I_j},n}(\textbf{x}_{I_j}^*)(1+o(1))$ in the last step.
This gives the expression for the variance in Lemma \ref{lemfd2}. For the calculation of the cumulants of $V^{-1/2}\hat{\theta}_{I_j}^{add,FD}$ we first note that the order of the variance $V=V(\hat{\theta}_{I_j}^{add,FD}(\textbf{x}_{I_j}))$ can be calculated in the same way as in the proof of Lemma \ref{lemfd}, which gives $V=O(n^{-d}h^{-d_j-2\beta_j}a_n^{-d_j})$. Therefore we have to show
\[ |cum_l(n^{d/2}h^{d_j/2+\beta}a_n^{d_j/2}\hat{\theta}_{I_j}^{add,FD})| = n^{ld/2}h^{l(d_j/2+\beta}a_n^{ld_j/2}|cum_l(\hat{\theta}_{I_j}^{add,FD})| \rightarrow 0 \]
for $l\ge 3$. By a straightforward calculation it follows that
\begin{eqnarray*}
&& |cum_l(\hat{\theta}_{I_j}^{FD})(\textbf{x}_{I_j}^*)|
\\ &=& \Big|\frac{1}{(2n+1)^{l(d-d_j)}n^{ld_j}h^{ld_j}a_n^{ld_j}}  \sum_{\textbf{k}_1,...,\textbf{k}_l \in \{-n,...,n\}^d}
 \prod_{m=1}^l \Big( \int_{\mathbb{R}^{d_j}}  e^{-i\langle \textbf{w},(\textbf{x}_{I_j}^*-z_{k_{m,I_j}})\rangle /h} \frac{\Phi_K(\textbf{w})}{\phi_\psi(\frac{\textbf{w}}{h})} d\textbf{w}\Big) cum(\epsilon_{\textbf{k}_1},...,\epsilon_{\textbf{k}_l}) \Big|
\\ &\le&\frac{C}{(2n+1)^{l(d-d_j)}n^{ld_j}h^{ld_j}a_n^{ld_j}}  \sum_{\textbf{k}_1,...,\textbf{k}_l \in \{-n,...,n\}^d}
 \prod_{m=1}^l  \Big( \int_{\mathbb{R}^{d_j}}   \frac{|\Phi_K(\textbf{w})|}{|\phi_\psi(\frac{\textbf{w}}{h})|} d\textbf{w}\Big)  |cum(\epsilon_{\textbf{k}_1},...,\epsilon_{\textbf{k}_l})|
\\ &=& \frac{C}{(2n+1)^{l(d-d_j)}n^{ld_j}h^{ld_j}a_n^{ld_j}} \frac{1}{h^{l\beta}} \sum_{\textbf{k}_1,...,\textbf{k}_l \in \{-n,...,n\}^d}
  |cum(\epsilon_{\textbf{k}_1},...,\epsilon_{\textbf{k}_l})|
\\ &=& \frac{C}{(2n+1)^{l(d-d_j)}n^{ld_j}h^{ld_j}a_n^{ld_j}} \frac{1}{h^{l\beta}}  (2n+1)^d,
\end{eqnarray*}
because by \eqref{MAq} $\textbf{k}_1$ can be chosen arbitrarily and $\textbf{k}_2,...,\textbf{k}_l$ have only $(4q+1)^d$ possibilities to be chosen and their bound is independent of $n$. Thus the $l$-th cumulant is of order $ n^{-ld/2+1}h^{-ld_j/2}a_n^{-ld_j/2},$ which converges to zero for $l \ge 3$. The result for $\hat{\theta}^{add,FD}$ follow immediately from the results of $\hat{\theta}_{I_j}^{FD}$.

\end{document}